\documentclass[12pt]{amsart}
\usepackage{amscd,amssymb,latexsym,hyperref}
\usepackage[graph,frame,poly,arc]{xy}  
\usepackage{subfigure}

\theoremstyle{plain}
\newtheorem{theorem}[subsection]{Theorem}
\newtheorem{lemma}[subsection]{Lemma}
\newtheorem{prop}[subsection]{Proposition}
\newtheorem{cor}[subsection]{Corollary}

\theoremstyle{definition}
\newtheorem{remark}[subsection]{Remark}
\newtheorem{definition}[subsection]{Definition}
\newtheorem{example}[subsection]{Example}

\numberwithin{equation}{section}

\newenvironment{alphenum}%
{%
   \renewcommand{\theenumi}{\alph{enumi}}%
   \renewcommand{\labelenumi}{(\theenumi)}%
   \begin{enumerate}
}
{%
   \end{enumerate}%
}

\newenvironment{romenum}%
{%
   \renewcommand{\theenumi}{\roman{enumi}}%
   \renewcommand{\labelenumi}{(\theenumi)}%
   \begin{enumerate}
}
{%
   \end{enumerate}%
}

{%
   \renewcommand{\theenumi}{\Alph{enumi}}%
   \renewcommand{\labelenumi}{(\theenumi)}%
   \begin{enumerate}
}%
{%
   \end{enumerate}%
}

\renewcommand{\arraystretch}{1.1}
\newcommand{\abs}[1]{\left\lvert #1\right\rvert}
\newcommand{\ind}[1]{\left\lvert #1\right\rvert}
\newcommand{\tsum}{\textstyle\sum\limits}

\newcommand{\triv}{\mathbf{1}}
\newcommand{\surj}{\twoheadrightarrow}
\newcommand{\inj}{\hookrightarrow}
\newcommand{\nor}{\triangleleft}

\newcommand{\A}{\mathcal{A}}

\newcommand{\LL}{\mathcal{L}}

\newcommand{\MM}{\mathcal{M}}

\newcommand{\C}{\mathbb{C}}
\newcommand{\F}{\mathbb{F}}
\newcommand{\Z}{\mathbb{Z}}
\newcommand{\Q}{\mathbb{Q}}
\newcommand{\R}{\mathbb{R}}
\newcommand{\RP}{\mathbb{RP}}
\newcommand{\K}{\mathbb{K}}
\newcommand{\SP}{\mathbb{S}}

\newcommand{\T}{\mathbb{T}}

\newcommand{\G}{\Gamma}
\newcommand{\s}{\sigma}
\renewcommand{\a}{{\alpha }}
\renewcommand{\b}{{\beta }}
\renewcommand{\l}{\lambda}
\renewcommand{\t}{\tau}
\renewcommand{\o}{\omega}
\renewcommand{\O}{\Omega}

\DeclareMathOperator{\rank}{rank}
\DeclareMathOperator{\corank}{corank}
\DeclareMathOperator{\coker}{coker}

\DeclareMathOperator{\Aut}{Aut}
\DeclareMathOperator{\Hom}{Hom}
\DeclareMathOperator{\End}{End}
\DeclareMathOperator{\Epi}{Epi}

\DeclareMathOperator{\GL}{GL}
\DeclareMathOperator{\Sym}{Sym}
\DeclareMathOperator{\Irrep}{Irrep}

\def\im{\operatorname{im}}
\DeclareMathOperator{\Tors}{Tors}

\DeclareMathOperator{\ii}{i}

\DeclareMathOperator{\ch}{char}
\DeclareMathOperator{\ord}{ord}
\newcommand{\Alt}{A}
\DeclareMathOperator{\Frat}{Frat}
\DeclareMathOperator{\partition}{\pi}
\DeclareMathOperator{\ann}{ann}

\newcommand{\sml}[1]{\text{\footnotesize{$#1$}}}

\begin{document}

\title[Hall invariants and characteristic varieties]%
{Hall invariants, homology of subgroups,\\ and characteristic varieties}

\author{Daniel~Matei}
\address{Department of Mathematics,
University of Rochester,
Rochester, NY 14627}
\email{\href{mailto:dmatei@math.rochester.edu}{dmatei@math.rochester.edu}}

\author{Alexander~I.~Suciu}
\address{Department of Mathematics,
Northeastern University,
Boston, MA 02115}
\email{\href{mailto:alexsuciu@neu.edu}{alexsuciu@neu.edu}}
\urladdr{\href{http://www.math.neu.edu/~suciu/}%
{http://www.math.neu.edu/\~{}suciu}}

\subjclass[2000]{Primary 20J05, 57M05; Secondary 20E07, 52C35}

\keywords{Hall invariant, Eulerian function, 
finite-index subgroup, homology, Alexander matrix, 
characteristic variety, torsion point, arrangement}


\thanks{Research partially supported by an RSDF grant from
Northeastern University.}

\begin{abstract}
Given a finitely-generated group $G$, and a finite group
$\Gamma$, Philip Hall defined $\delta_\Gamma(G)$ to be the
number of factor groups of $G$ that are isomorphic to $\Gamma$.
We show how to compute the Hall invariants
by cohomological and combinatorial methods, when $G$ is
finitely-presented, and $\Gamma$ belongs to a certain class of
metabelian groups.  Key to this approach is the stratification of
the character variety, $\Hom(G,\K^*)$, by the jumping loci of the
cohomology of $G$, with coefficients in rank~$1$ local systems
over a suitably chosen field $\K$.  Counting relevant torsion
points on these ``characteristic" subvarieties gives
$\delta_\Gamma(G)$.  In the process, we compute the 
distribution of prime-index, normal subgroups $K\nor G$ 
according to $\dim_{\K} H_1(K;\K)$, provided 
$\ch \K\ne \ind{G:K}$. In turn, we use this distribution 
to count low-index subgroups of $G$.  We illustrate these 
techniques in the case when $G$ is the fundamental group of 
the complement of an arrangement of either affine lines in 
$\C^{2}$, or transverse planes in $\R^4$.
\end{abstract}

\maketitle

\tableofcontents

\section{Introduction}
\label{sec:intro}

\subsection{Hall invariants}
In \cite{HaP}, Philip Hall introduced several 
notions in group theory.   Given a finitely-generated group $G$, 
and a finite group $\Gamma$, he defined $\delta_{\Gamma}(G)$
to be the number of surjective representations of $G$ to $\Gamma$,
up to automorphisms of $\Gamma$:
\begin{equation}
\label{eq:delgamma}
\delta_{\Gamma}(G)=\abs{\Epi(G,\Gamma)/\Aut\Gamma}.
\end{equation}
In other words, the {\em Hall invariant} $\delta_{\Gamma}(G)$ 
counts all factor groups of $G$ that are isomorphic to $\Gamma$. 
When $G=\pi_1(X)$ is the fundamental group of a connected $2$-complex
$X$ with finite $1$-skeleton, $\delta_{\Gamma}(\pi_1(X))$ counts all 
connected, regular covers of $X$  with deck transformation group $\G$.

Suppose $G$ has a finite presentation, with generators $x_1,\dots ,x_n$, 
and relators $r_1,\dots ,r_m$.  The Hall invariant $\delta_{\Gamma}(G)$ can 
be computed by counting generating sets of $\G$ that have size $n$ and satisfy 
the relations $r_j$ in $\G$, and then dividing the result by the order of
$\Aut(\Gamma)$.   While this method can be implemented on a computer 
algebra system like {\sl GAP} \cite{gap}\footnote{%
The {\sl GAP}  command that computes $\delta_{\G}(G)$ is 
\textsf{GQuotients(G,$\G$)}.  We were led to consider the Hall 
invariants after reading about this command in the {\sl GAP}
manual.},   the  computation breaks down even for moderately large $n$
and $\abs{\G}$.   One of the purposes of this paper is to show how the Hall
invariants $\delta_{\Gamma}(G)$ can be computed more efficiently, by combinatorial
and homological methods, at least when $\Gamma$ belongs to a certain  
class of finite metabelian groups.

\subsection{Abelian representations}
\label{subsec:abel}
In Section~\ref{sec:euler}, we start by recalling a formula of 
Hall~\cite{HaP}.  Using what he called the {\em Eulerian function} of the 
finite group $\G$, together with  M\"{o}bius inversion, Hall showed that 
$\delta_{\Gamma}(G)=\tfrac{1}{\abs{\Aut\G}}\sum_{H\le\G}\mu(H)\abs{\Hom(G,H)}$, 
where $\mu$ is the M\"{o}bius function of the subgroup lattice of $\G$.  
For a $p$-group $\G$, the M\"{o}bius and Eulerian functions were computed by
Weisner~\cite{We}.    

In Section~\ref{sec:abelian}, we use these results of Hall and Weisner, 
together with a result from Macdonald's book \cite{Mac}, 
to arrive at a completely explicit formula for $\delta_{\Gamma}(G)$ 
in the simplest case:  that of a finite abelian group $\G$. 
The expression for $\delta_{\Gamma}(G)$, given in 
Theorem~\ref{thm:abelcount}, depends only on $\G$ 
and the abelian factors of $G$.  A similar expression  
was obtained (by other means) in \cite{KCL}, in 
the particular case $G=F_n$.  

\subsection{Metabelian representations}
\label{subsec:metabel}
The next level of difficulty in computing $\G$-Hall invariants 
is presented by (split) metabelian groups $\G$.  
Suppose $\G=B\rtimes_{\sigma}C$ is a semidirect product 
of abelian  groups, with monodromy homomorphism $\s:C\to\Aut(B)$.    
An epimorphism $\l:G\surj\G$ may be thought of as the lift of an epimorphism
$\rho:G\surj C$. As explained in~\cite{Liv}, the lifts of a fixed homomorphism 
$\rho:G\to C$ are parametrized by $H^1(G,B_{\rho})$, the first cohomology group of
$G$  with coefficients in the $G$-module $B=B_{\rho}$, with action 
$g\cdot b=\sigma(\rho(g))(b)$. 

We consider split metabelian groups $\G=B\rtimes_{\sigma}C$ 
for which $B$ is an elementary abelian group, and $C$ is a finite cyclic group. 
If $\K$ is a finite field with additive group $B$, then
$H^1(G,B_{\rho})$ may be identified with $H^1(G,\K_{\rho})$ where again
$G$ acts on $\K$ by means of $\rho$. Now Shapiro's Lemma identifies  
the twisted cohomology group $H^1(G,\K_{\rho})$ with 
the untwisted cohomology group $H^1(K_{\rho},\K)$, 
where $K_{\rho}=\ker\rho$. We are thus led to investigate 
the homology of finite-index, normal subgroups of $G$.

\subsection{Homology of finite-index subgroups}
\label{subsec:homfinind}
Let $G$ be a finitely-presented group, and $K\nor G$ a normal 
subgroup.  Assume the quotient group, $\G=G/K$, is finite. 
A procedure to compute $H_1(K,\Z)$ from a presentation of $G$ 
and the coset representation of $G$ on $G/K$ was given by Fox~\cite{Fx1}. 
The efficiency of Fox's method decreases rapidly with the
increase in the index $\abs{\G}=\ind{G:K}$. In Section~\ref{sec:homology}, 
we overcome this problem, at least partially.  Our approach 
(similar to that of Hempel~\cite{He1, He2} and Sakuma~\cite{Sa}) 
is based on the representation theory of $\G$, over suitably 
chosen fields $\K$.     

Consider the homology group $H_1(K,\K)$, and set 
$b_1^{(q)}(K):=\dim_{\K} H_1(K,\K)$, where $q=\ch\K$. 
The idea is to break $H_1(K,\K)$ into a direct sum, according to the 
decomposition of the group algebra $\K\G$ into irreducible 
representations. In order for this to work, we need the field $\K$ 
to be ``sufficiently large" with respect to $\G$; that is,  
$q$ should not to divide $\abs{\G}$, and $\K$ should contain 
all roots of unity of order equal to the exponent of $\G$.  
Let $\lambda:G\surj\Gamma$ be an epimorphism, with kernel $K=K_{\l}$.
In Theorem~\ref{thm:genlibsak}, we prove:
\begin{equation}
\label{eq:b1decomp}
b_1^{(q)}(K_{\l}) = b_1^{(q)}(G)+\sum_{\rho\ne \triv}
n_{\rho}( \corank J^{\rho\circ \lambda}- n_{\rho}),
\end{equation}
where $\rho$ runs through all non-trivial, irreducible 
representations of $\Gamma$ over the field $\K$, and 
$J^{\rho\circ \lambda}$ is the Jacobian matrix of Fox derivatives 
of the relators, $J=J_G$, followed by the representation 
$\rho\circ \lambda: G\to \GL(n_{\rho},\K)$.

When $\G$ is abelian and $\K=\C$, we 
recover from \eqref{eq:b1decomp} a well-known result 
of Libgober~\cite{Li1} and Sakuma~\cite{Sa}: 
$b_1(K_{\l}) = b_1(G) +\sum_{\rho\ne\triv}
( \corank J^{\rho\circ \lambda}- 1)$, 
where $\rho$ runs through all non-trivial, irreducible,
complex representations of $\Gamma$.  For other choices 
of $\K$, formula \eqref{eq:b1decomp} gives information 
about the $q$-torsion coefficients of $H_1(K,\Z)$, 
provided $q\nmid \abs{\G}$.

\subsection{Torsion points on characteristic varieties}
\label{subsec:chvars}
The next step is to interpret formula \eqref{eq:b1decomp} in terms of 
the ``Alexander stratification" of the character variety of $G$. 
This can be done for an arbitrary finite abelian group $\G$, 
but, for simplicity, we restrict our attention to the case 
when $\G$ is cyclic, which is enough for our purposes here.

In Section \ref{sec:alexcv}, we start by reviewing the pertinent 
material on Alexander ideals and their associated varieties, 
in a more general context than usual.  
The Alexander matrix, $A_G$, is the abelianization of  
the Fox Jacobian, $J_G$.   The {\em $d$-th characteristic
variety}, $V_d(G,\K)$, is the subvariety of $\Hom(G,\K^*)$ 
defined by the  codimension~$d$ minors of $A_G$. 
It can be shown that $V_d(G,\K)\setminus \{\mathbf{1}\}$ 
is the set of non-trivial characters $\mathbf{t}\in\Hom(G,\K^*)$  
for which $\dim_{\K} H^1(G,\K_{\mathbf{t}}) \ge d$, see \cite{Hr} and
Remark~\ref{rem:hironaka}. 

In Section~\ref{sec:bettitors}, we study the relationship  
between torsion points on the characteristic varieties of $G$ 
and the homology of finite-index, normal subgroups $K\nor G$. 
As mentioned above, we only consider the case when $\G=G/K$ 
is cyclic, say $\G=\Z_N$.  In Theorem~\ref{thm:homcyclic}, we prove:
\begin{equation}
\label{eq:b1jp}
b_1^{(q)}(K_{\l}) = b_1^{(q)}(G)+\sum_{1\ne k\mid N} \phi(k) d_{\K}(\l^{N/k}).
\end{equation}
where $d_{\K}(\mathbf{t})=\max\, \{ d \mid \mathbf{t}\in V_d(G,\K)\}$ 
is the {\em depth} of the character $\mathbf{t}\in \Hom(G,\K^*)$ 
with respect to the Alexander stratification.  
In particular, if $N=p$ is prime, then: 
\begin{equation}
\label{eq:primejump}
b_1^{(q)}(K_{\l}) = b_1^{(q)}(G)+(p-1)d_{\K}(\l).
\end{equation}

In view of \eqref{eq:primejump}, we define $\b_{p,d}^{(q)}(G)$ 
to be the number of index~$p$, normal subgroups $K\nor G$ 
for which $ b_1^{(q)}(K) = b_1^{(q)}(G) + (p-1) d$. 
In Theorem~\ref{thm:torscount}, we prove:
\begin{equation}
\label{eq:betapdq}
\b_{p,d}^{(q)}(G)=  \tfrac{1}{p-1}\,\abs{\Tors_{p,d}(G,\K)
\setminus \Tors_{p,d+1}(G,\K)},
\end{equation}
where $\K=\C$ if $q=0$, and $\K=\F_{q^s}$
($s=\text{order of $q$ in $\F_p^*$}$) if $q\ne 0$, and 
$\Tors_{p,d}(G,\K)$ is
the set of characters in $V_d(G,\K)$ of order exactly $p$. 

\subsection{Metabelian Hall invariants and low-index subgroups}
\label{subsec:metabelhall}
Once this is done, we are ready to return to the Hall invariants of $G$.  
In Section~\ref{sec:metacyclic}, we compute $\delta_{\G}(G)$, for 
split metabelian groups $\G$ of the form 
$M_{p,q^s}=\Z_q^s\rtimes_{\s} \Z_p$, where $p$ and $q$ are distinct primes,
$s=\ord_p(q)$, and $\s$ has order exactly~$p$.      
Examples are the dihedral groups $D_{2p}=M_{2,p}$ and the alternating group
$A_4=M_{3,4}$.   In Theorem~\ref{thm:delmpqs}, we prove:
\begin{equation}
\label{eq:delcount}
\delta_{M_{p,q^s}}(G) =\tfrac{p-1}{s(q^s-1)} \sum_{d\ge 1}
\b_{p,d}^{(q)}(G) (q^{sd}-1).
\end{equation}
This generalizes a result of Fox \cite{Fx2}, who was the first to use 
Alexander matrices for counting metacyclic representations 
of fundamental groups of knots and links.  Put together, formulas
\eqref{eq:delcount} and \eqref{eq:betapdq} express the Hall invariant 
$\delta_{M_{p,q^s}}(G)$ in terms of the number of $p$-torsion points 
on the Alexander strata of the character variety $\Hom(G,\F_{q^s})$.  

In Section~\ref{sec:lis}, we use formula \eqref{eq:delcount}, 
together with several formulas from \S\S\ref{sec:euler}--\ref{sec:abelian}, 
to derive information about the number, $a_k(G)$, 
of index $k$ subgroups of $G$.  It was Marshall Hall \cite{Ha} 
who showed how to compute these numbers recursively, in terms of 
$\abs{\Hom(G,S_l)}$, $1\le l\le k$.  
Applying this method, we obtain (in Theorem~\ref{thm:index2and3}):
\begin{equation}
\label{eq:a3g}
a_3(G)=\tfrac{1}{2}(3^{n}-1) + \tfrac{3}{2} 
\sum_{d\ge 1} \beta_{2,d}^{(3)}(G) (3^d-1),
\end{equation} 
where $n=b_1^{(3)}(G)$.  
We also give formulas of this sort for the number,
$a_k^{\nor}(G)$,  of index $k$, normal subgroups of $G$, 
provided $k\le 15$ and $k\ne 8$ or $12$.

\subsection{Arrangement groups}
\label{subsec:arrgp}
We conclude with some explicit examples and computations 
in the case when $G$ is the fundamental group of the complement 
of a subspace arrangement.  This is meant to illustrate the theory 
developed so far, in a setting where topology and combinatorics 
are closely intertwined. 

In Section~\ref{sec:complexarr}, we look at complex hyperplane 
arrangements.  By the Lefschetz-type theorem of Hamm and L\^{e}, 
it is enough to consider arrangements of affine lines in $\C^2$.  
If $G$ is the group of such an arrangement, 
the characteristic varieties $V_d(G,\C)$ are well
understood:  they consist of subtori of the character torus, possibly 
translated by roots of unity.   Furthermore, the tangent cone 
at the origin to $V_d(G,\C)$ coincides with the ``resonance" 
variety $R_d(G,\C)$, which is determined by the combinatorics 
of the arrangement.  The components of $V_d(G,\C)$ not 
passing  through the origin, though, are not {\it a priori} 
combinatorially determined.  Their appearance affects 
the torsion coefficients in the homology of certain finite abelian covers 
of the complement, as we show in Example~\ref{ex:a22}.

In Section~\ref{sec:planes}, we turn to real arrangements.  
More precisely, we consider arrangements of 
transverse planes through the  origin of $\R^4$.  If $G$ 
is the group of such a {\em $2$-arrangement}, the varieties $V_d(G,\C)$ 
need not be unions of translated subtori, as shown in \cite{MS1}, 
and also here, in Example~\ref{ex:a31425}.   Furthermore, 
the tangent cone at the origin to $V_d(G,\C)$ may not coincide 
with the resonance variety $R_d(G,\C)$, as we  point out in
Remark~\ref{rem:tcone}.  Finally, using the metabelian  Hall 
invariants $\delta_{S_3}$ and $\delta_{\Alt_4}$, we  recover the 
homotopy-type classification of complements  of $2$-arrangements 
of $n\le 6$ planes in $\R^4$ (first  established in \cite{MS1}), 
and extend it to horizontal arrangements of $n=7$ planes.  

\subsection*{Acknowledgment}  
The computations for this work were done with the help of the
packages {\sl GAP~4.1} \cite{gap},
{\sl Macaulay~2} \cite{GS}, and {\sl  Mathematica~4.0}.

\section{Eulerian functions and Hall invariants}
\label{sec:euler}

We start by reviewing two basic notions introduced by Philip Hall
in~\cite{HaP}:   the Eulerian function, $\phi(\G,n)$, of a finite group
$\G$, and the Hall invariants, $\delta_{\Gamma}(G)$,
of a finitely-generated group $G$.

\subsection{Eulerian function}
\label{subsec:euler}
Let $\Gamma$ be a finite group.
The {\em Eulerian function} of $\G$ is defined as
\begin{equation}
\label{eq:euler}
\phi(\G,n)=\#\{\text{ordered $n$-tuples $(g_1,\dots ,g_n)$ that 
generate $\G$}\},
\end{equation}
where repetitions among the $g_i$'s are allowed.  For example,
$\phi(\Z_k,1)=\phi(k)$, the usual Euler totient function.

Let $L(\G)$ be the lattice of subgroups of $\G$, ordered by inclusion.
Let $\mu: L(\G)\to \Z$ be the {\em M\"{o}bius function}, defined 
inductively by $\mu(\G)=1$, $\sum_{H\le K} \mu(H)=0$.
Then, the Eulerian function of $\Gamma$ is given by:
\begin{equation}
\label{eq:phall}
\phi(\G,n)=\sum_{H\le \G} \mu(H) \abs{H}^n,
\end{equation}
see \cite{HaP}, and also \cite{Br} for a recent account.

The Eulerian function and the M\"{o}bius function of 
a finite $p$-group were computed by Weisner in~\cite{We}. 
To formulate Weisner's results, recall that the 
{\em Frattini subgroup} of a finite group $\G$, 
denoted $\Frat\G$, is the intersection of all maximal, 
proper subgroups of $\G$.   
If $\G$ is a $p$-group, then $\Frat\G=[\G,\G]\cdot \G^p$, 
by the Burnside Basis Theorem (cf.~\cite[p.~140]{Rob}). 

According to Weisner, the M\"{o}bius function of a 
finite $p$-group $\G$ is given by:  
\begin{equation}
\label{eq:mobiuspg}
\mu(H)=
\begin{cases}
(-1)^{d} p^{\frac{d(d-1)}{2}}, \text{ where } 
p^{d}=\ind{\G:H} &\text{ if } \Frat\G\le H,\\
0  &\text{ if } \Frat\G\not\le H.
\end{cases}
\end{equation}
Now set $p^{r}=\abs{\G}$ and $p^{s}=\ind{\G:\Frat\G}$.
The Eulerian function of $\G$ is then given by:
\begin{equation}
\label{eq:phipg}
\phi(\G,n)=p^{(r-s) n} \prod_{i=0}^{s-1} (p^{n}-p^{i}).
\end{equation}

\subsection{Hall invariants}
\label{subsec:hallinv}
Let $G$ be a finitely-generated group, and $\G$ a finite group.  
Let $\sigma_{\G}(G)=\abs{\Hom(G,\Gamma)}$ be the number 
of homomorphisms $G\to\G$, and 
$\phi_{\G}(G)=\abs{\Epi(G,\Gamma)}$ the number 
of epimorphisms $G\surj\G$. 
The relation between $\sigma$ and $\phi$ is 
given by Hall's enumeration principle:
\begin{equation}
\label{eq:hallenum1}
\sigma_{\G}(G)=\sum_{H\le\G}\phi_{H}(G),
\end{equation}
or, by M\"obius inversion:
\begin{equation}
\label{eq:hallenum2}
\phi_{\G}(G)=\sum_{H\le\G}\mu(H)\sigma_{H}(G).
\end{equation}

\begin{definition}
\label{def:hallinv}
Let $G$ be a finitely-generated group. Let $\G$ be a finite group, 
with automorphism group $\Aut \G$. 
The {\em $\G$-Hall invariant} of $G$ is 
$\delta_{\Gamma}(G)=\phi_{\G}(G)/\abs{\Aut\Gamma}$.  
\end{definition}

Since $\Aut\Gamma$ acts freely and transitively on $\Epi(G,\Gamma)$,
the number $\delta_{\Gamma}(G)$ is an integer, which counts 
epimorphisms $G\surj\Gamma$, up to automorphisms of $\Gamma$.  
In other words,\ $\delta_{\Gamma}(G)$ is the number of homomorphs 
of $G$ that are isomorphic to $\G$.  

Note that
\begin{equation}
\label{eq:product}
\phi_{\Gamma_1\times \Gamma_2}(G)=\phi_{\Gamma_1}(G) \phi_{\Gamma_2}(G)
\end{equation}
whenever $\G_1$ and $\G_2$ are finite groups, with
$(\abs{\G_1},\abs{\G_2})=1$.  In that situation, we 
also have $\Aut(\Gamma_1\times \Gamma_2)=\Aut(\Gamma_1)\times\Aut(\Gamma_2)$, and so 
$\delta_{\Gamma_1\times \Gamma_2}(G)=\delta_{\Gamma_1}(G)
\delta_{\Gamma_2}(G)$.

Now let $G=F_n$, the free group of rank $n$. Clearly, 
$\sigma_{\G}(F_n)=\abs{\G}^n$ and  
$\phi_{\G}(F_n)=\phi(\G,n)$.  Hence, by \eqref{eq:phall}:
\begin{equation}
\label{eq:delfree}
\delta_{\G}(F_n)=\frac{\sum_{H\le \G} \mu(H) \abs{H}^n}{\abs{\Aut(\G)}}.
\end{equation}

\begin{figure}
\setlength{\unitlength}{1cm}
\subfigure{%
\label{fig:mobmpq}%
\begin{minipage}[t]{0.38\textwidth}
\begin{picture}(4,4.8)(1,-3.2)
\xygraph{!{0;<9.5mm,0mm>:<0mm,10.5mm>::}  
[]*L(7.5){\sml{1}}*+=[o]+[F]{\text{\small\!$M_{p,q}$\!}}
(-[dlll]*R(3){\sml{-1}}*+=[o]+[F]{\,\Z_q\,}
(-[dddrrr]*R(5){\sml{q}}*+=[o]+[F]{1})
,-[ddl]*DL(2.5){\!\sml{-1}}*+=[o]+[F]{\,\Z_p\,}
(-[ddr]{\phantom{1}})
,-[dd]*DL(2.5){\sml{-1}}%
*L(2.0){\quad \ldots}%
*U(2.0){\quad\ \ \:\underbrace{\qquad\qquad\quad\ }_{q}}%
*+=[o]+[F]{\,\Z_p\,}
(-[dd]{\phantom{1}})
,-[ddrr]*DL(2.5){\sml{-1}}*+=[o]+[F]{\,\Z_p\,}
(-[ddll]{\phantom{1}})
) 
}
\end{picture}
\end{minipage}
}
\setlength{\unitlength}{0.8cm}
\subfigure{%
\label{fig:moba4}%
\begin{minipage}[t]{0.38\textwidth}
\begin{picture}(4,4.8)(0,-4.2)
\xygraph{!{0;<9.6mm,0mm>:<0mm,8.5mm>::}
[]*L(6.5){\sml{1}}*+=[o]+[F]{A_4}
(-[dl]*R(3){\sml{-1}}*+=[o]+[F]{\,\Z_2^2\,}
(
-[ddl]*UR(3.35){\sml{0}}*+=[o]+[F]{\,\Z_2\,}
(
-[ddrr]*R(5){\sml{4}}*+=[o]+[F]{1}
)
,-[dd]*UR(3.35){\sml{0}}*+=[o]+[F]{\,\Z_2\,}
(-[ddr]{\phantom{1}})
,-[ddr]*UR(3.35){\sml{0}}*+=[o]+[F]{\,\Z_2\,}
(-[dd]{\phantom{1}})
)
,-[ddr]*U(3.75){\sml{\, -1}}*+=[o]+[F]{\,\Z_3\,}
(-[dddl]{\phantom{1}})
,-[ddrr]*U(3.75){\sml{\! -1}}*+=[o]+[F]{\,\Z_3\,}
(-[dddll]{\phantom{1}})
,-[ddrrr]*U(3.75){\sml{\!\! -1}\quad }*+=[o]+[F]{\,\Z_3\,}
(-[dddlll]{\phantom{1}})
,-[ddrrrr]*U(3.75){\:\sml{-1}}*+=[o]+[F]{\,\Z_3\,}
(-[dddllll]{\phantom{1}})
) 
}
\end{picture}
\end{minipage}
}
\caption{\textsf{The subgroup lattice and 
M\"obius function of $M_{p,q}$ and $\Alt_4$}}
\label{fig:sgplattice}
\end{figure}

\begin{example}
\label{ex:deltametafree}
Let $M_{p,q}=\Z_q\rtimes \Z_p$ be the metacyclic group of order $pq$,
where $p$ and $q$ are primes, with $p\mid (q-1)$.
Its subgroup lattice and M\"{o}bius function are 
shown in Figure~\ref{fig:sgplattice}. 
The automorphism group of $M_{p,q}$ is isomorphic to 
$M_{q-1,q}\cong\Z_q\rtimes \Aut(\Z_q)$, 
the holomorph of $\Z_q$.  By \eqref{eq:delfree}:
\begin{equation}
\label{eq:deltameta}
\delta_{M_{p,q}}(F_n) = \frac{(p^{n}-1)(q^{n-1}-1)}{q-1}.
\end{equation}
\end{example}

\begin{example}
\label{ex:deltaltfree}
Let $\Alt_4$ be the alternating group on $4$ symbols.  
The M\"{o}bius function is given in Figure~\ref{fig:sgplattice}.  
Furthermore,  $\Aut(\Alt_4)\cong S_4$, the symmetric group on $4$ symbols.
We get:
\begin{equation}
\label{eq:deltalt}
\delta_{\Alt_4}(F_n) = \frac{(3^{n}-1)(4^{n-1}-1)}{6}.
\end{equation}
\end{example}

\section{Counting abelian representations}
\label{sec:abelian}

In this section, we show how to compute the Hall
invariant $\delta_{\Gamma}(G)$, in case $\G$ is a finite abelian
group.  We start with the well-known computation of the order of $\Aut(\G)$. 

For a prime $p$, denote by $\G_p$ the $p$-torsion part of $\G$.  
Then $\G=\bigoplus_{p\mid \abs{\G}} \G_p$, and 
\begin{equation}
\label{eq:prodaut}
\abs{\Aut(\G)}=\prod\limits_{p \mid \, \abs{\G}} \abs{\Aut(\G_p)}.
\end{equation}

Let $A$ be a (finite) abelian $p$-group.  Then 
$A=\Z_{p^{\pi_1}}\oplus\cdots\oplus\Z_{p^{\pi_r}}$, for some 
positive integers $\pi_1\ge \dots \ge \pi_r$, and so 
$A$ determines (and is determined by) 
a partition $\partition(A)=(\pi_1, \dots, \pi_r)$. 
Given such a partition $\pi$, let $l(\pi)=r$ be its length, 
and $|\pi|=\sum_{i=1}^{r} \pi_i$ its weight.  Also, let 
$\langle\pi\rangle=\sum_{i=1}^{r}(i-1)\pi_i$. 
 Then:
\begin{equation}
\label{eq:autpgp}
\abs{\Aut(\Z_{p^{\pi_1}}\oplus \cdots\oplus\Z_{p^{\pi_r}})}
=p^{|\pi| +2\langle \pi \rangle}\prod_{k\ge 1}\varphi_{m_k(\pi)}
(p^{-1}),
\end{equation}
where 
$m_k(\pi)=\#\{ j\mid \pi_j=k\}$ is the multiplicity of $k$ in $\l$, 
and $\varphi_m(t)=\prod_{i=1}^{m} (1-t^i)$; 
see Macdonald~\cite[p.~181]{Mac}\footnote{We are 
grateful to A.~Zelevinsky for pointing this reference to us.}.

Given a partition $\l$, let $\l'$ be the partition with $\l'_i=\l_i-1$.
If $\tau$ is another partition, let
$\theta_i(\l,\tau)=\sum_{j=1}^{l(\tau)}\min(\l_i,\tau_{j})$, for 
$1\le i \le l(\l)$, and
$\theta(\l,\tau)=\sum_{i=1}^{l(\l)}\theta_i(\l,\tau)$.
With these notations, we have the following:

\begin{theorem}
\label{thm:abelcount}
Let $G$ be a finitely-generated group  and $\G$ a
finite abelian group.   Write $H_1(G)=\Z^n\oplus T$, where 
$T$ is finite. For each prime $p$ dividing $\abs{\G}$, 
let $\l=\partition(\G_p)$ and $\tau=\partition(T_p)$ 
be the corresponding partitions.  Then:
\begin{equation}
\label{eq:deltapgroup}
\delta_{\G}(G)=\prod_{p \mid \, \abs{\G}} \frac
{p^{(|\l|-l(\l)) n + \theta(\l',\tau)}
\prod_{i=1}^{l(\l)}
(p^{n+\theta_i(\l,\tau)-\theta_i(\l',\tau)}-p^{i-1})}
{p^{|\l|+2\langle \l \rangle}\prod_{k\ge 1}\varphi_{m_k(\l)}
(p^{-1})}.
\end{equation}
\end{theorem}

\begin{proof}
Since $\G$ is abelian, every homomorphism $G\to\G$
factors through $H_1(G)=\Z^n\oplus T$.  By \eqref{eq:product}, 
we have $\phi_{\Gamma}(G)=\prod_{p \mid \, \abs{\G}}
\phi_{\Gamma_p}(\Z^n\oplus T)$.  Furthermore, every homomorphism 
$\Z^n\oplus T \to \G_p$ factors through $\Z^n\oplus T_p$. Hence:
\begin{equation}
\label{eq:prodabel}
\phi_{\Gamma}(G)=\prod_{p \mid \, \abs{\G}}
\phi_{\Gamma_p}(\Z^n\oplus T_p).
\end{equation}
By Hall's enumeration principle \eqref{eq:hallenum2}, we have:
\begin{equation}
\label{eq:phigamma}
\phi_{\Gamma_p}(\Z^n\oplus
T_p)=\sum_{H\le\G_p}\mu(H)\sigma_{H}(\Z^n\oplus T_p).
\end{equation}   
Clearly, $\sigma_{H}(\Z^n \oplus T_p)=\abs{H}^n p^{\theta(\nu,\tau)}$, 
where $\nu=\pi(H)$. The M\"{o}bius function of a subgroup $H$ 
of $\G_p$ can be computed from Weisner's formula \eqref{eq:mobiuspg}, 
as follows.

Let $\Frat\G_p$ be the Frattini subgroup of $\G_p$.  Then 
$\Frat\G_p=(\G_p)^p$, and so the associated partition is 
$\l'$, where $\l=\pi(\G_p)$ and $\l'_i=\l_i-1$. 
If a subgroup $H$ of $\G_p$ 
does not contain $\Frat\G_p$, then $\mu(H)=0$.  
If $H$  contains $\Frat\G_p$, then the 
partition $\nu=\pi(H)$ is between $\l'$ and $\l$, 
i.e., $\nu_j=\l_j-1$ or $\l_j$.   
Order the set $\{j\mid \nu_j=\l_j-1\}$ as $(i_1\ge \cdots \ge i_d)$. 
Clearly, $p^d=\ind{\G:H}$, and so $\mu(H)=(-1)^{d} p^{d(d-1)/2}$. 

The number of subgroups $H\le \G_p$ such that
$\Frat\G_p \le H$ and $\pi(H)=\nu$ equals $p^a$, where 
$a=\sum_{k=1}^{d} (i_k-k)$. Set $d_{\nu}=d$ and  $a_{\nu}=a$.  
A simple calculation now shows:
\begin{align*}
\phi_{\G_p}(\Z^n \oplus T_p) 
&=
\sum_{\Frat\G_p \le H\le\G_p}\mu(H)\sigma_{H}(\Z^n \oplus T_p)\\
&=
\sum_{\l'\preceq\nu\preceq\l} (-1)^{d_{\nu}}
p^{d_{\nu}(d_{\nu}-1)/2}
p^{(|\l|-d_{\nu})n}
p^{\theta(\nu,\tau)} p^{a_{\nu}}\\
&=
p^{(|\l|-l(\l)) n + \theta(\l',\tau)}
\prod_{i\ge 1}
(p^{n+\theta_i(\l,\tau)-\theta_i(\l',\tau)}-p^{i-1}).
\end{align*}
This, together with formulas \eqref{eq:prodaut}, \eqref{eq:autpgp}, and 
\eqref{eq:prodabel}, yields \eqref{eq:deltapgroup}.
\end{proof}

Especially simple is the case when the group $G$ has 
torsion-free abelianization.

\begin{cor}
\label{cor:abelcountfree}
Let $G$ be a finitely generated group with
$H_1(G)=\Z^n$, and $\G$ a finite abelian
group. Write $\G=\prod\limits_{p \mid \, \abs{\G}} \G_p$,
where
$\G_p=\Z_{p^{\lambda_1}}\oplus\cdots\oplus\Z_{p^{\lambda_r}}$.
Then
\begin{equation}
\label{eq:delpgroup}
\delta_{\G}(G)=\prod_{p \mid \, \abs{\G}}
\frac{p^{|\l|(n-1)-2\langle\l\rangle} \varphi_n(p^{-1})}
{\varphi_{n-r}(p^{-1})\prod_{k\ge 1} \varphi_{m_k(\l)} (p^{-1})}.
\end{equation}
\end{cor}

A formula similar to \eqref{eq:delpgroup} was obtained by 
Kwak, Chun, and Lee~\cite[Theorem 3.4]{KCL}.  In particular 
(writing $\Z_p^{s}$ for the direct sum of $s$ copies of $\Z_p$):
\[
\delta_{\Z_{p^{s}}}(G)=\tfrac{p^{sn}-p^{(s-1)n}}{p^{s}-p^{s-1}}
,\quad
\delta_{\Z_{p}^{s}}(G)=\prod_{i=0}^{s-1}\tfrac{p^{n}-p^{i}}
{p^{s}-p^{i}}
,\quad
\delta_{\Z_{p}\oplus \Z_{p^{s}}}(G)=
\tfrac{(p^{s n}-p^{(s-1)
n})(p^{n}-p)}{p^{s+1}(p-1)^{2}}
.
\]


\section{Homology of finite-index subgroups}
\label{sec:homology}

In this section, we give a formula for computing the
first homology (with coefficients in a ``sufficiently large"
field) of finite-index, normal subgroups of a
finitely-presented group.

\subsection{Fox calculus}
\label{subsec:fox}

Let $G=\langle x_1,\dots ,x_{\ell}\mid r_1,\dots ,r_m\rangle$
be a finite presentation for the group $G$.
Let $F_{\ell}$ be the free group with
generators $x_1,\dots,x_{\ell}$, and $\phi:F_{\ell}\to G$
the presenting epimorphism.  Let $\Z F_{\ell}$ be the
group-ring of $F_{\ell}$, and $\epsilon:\Z F_{\ell}\to \Z$
the augmentation map.  For each $1\le j \le \ell$, there is a
Fox derivative, $\frac{\partial}{\partial x_j}: \Z F_{\ell}\to \Z F_{\ell}$,
which is the linear operator defined by the rules
$\frac{\partial 1}{\partial x_j}=0$,
$\frac{\partial x_i}{\partial x_j}=\delta_{ij}$, and
$\frac{\partial(uv)}{\partial x_j}=\frac{\partial u}{\partial x_j}\epsilon(v)+
u\frac{\partial v}{\partial x_j}$.  

Let $X$ be the $2$-complex associated with the presentation
$G=\langle x_1,\dots , x_{\ell} \mid r_1,\dots ,r_m\rangle$.
Let $\widetilde{X}$ be the universal cover,
and $C_*(\widetilde{X})$  
its augmented cellular chain complex. Picking as
generators for the chain groups the lifts of the cells
of $X$, the complex $C_*(\widetilde{X})$ becomes identified with
\begin{equation}
\label{eq:chaincomplex}
(\Z G)^m \xrightarrow{J_G} (\Z G)^{\ell}
\xrightarrow{\partial_1}\Z G\xrightarrow{\epsilon}\Z\to 0,
\end{equation}
where 
$\partial_1=
\begin{pmatrix}x_1-1& \cdots &x_{\ell}-1
\end{pmatrix}^{{\scriptscriptstyle{\top}}}$,
and \[J_G=\begin{pmatrix}
\frac{\partial r_i}{\partial x_j}\end{pmatrix}^{\phi}\]
is the Jacobian matrix of $G$, obtained by applying
the linear extension $\phi:\Z F_{\ell}\to \Z G$
to the Fox derivatives of the relators.

Clearly, the integral $m\times \ell$ matrix $J_G^{\epsilon}$ 
is a presentation matrix for $H_1(G)$. More generally, 
the abelianization of a finite-index subgroup $K\le G$ 
is given by the following result of Fox.

\begin{theorem}[Fox~\cite{Fx1}]
\label{thm:fox}
Let $G$ be a finitely-presented group, and $K<G$ a subgroup 
of index~$k$.  Let $J=J_G$ be the Jacobian matrix, 
$\sigma:G\to \Sym(G/K)\cong S_k$ the coset representation, 
and $\pi:S_k\to \GL(k,\Z)$ the permutation representation. 
Then $J^{\pi\circ \sigma}$ is a presentation 
matrix for the abelian group $H_1(K)\oplus\Z^{k-1}$.
\end{theorem}

\begin{proof}
By Shapiro's Lemma (cf.~\cite{Br82}), $H_1(K,\Z)$ is isomorphic to
$H_1(G,\Z[G/K])$, the first homology of the chain 
complex \eqref{eq:chaincomplex},
tensored over $\Z G$ with the module $\Z[G/K]$.
Under the identification $\Z [G/K]=\Z^k$, 
the boundary map $J\otimes \Z [G/K]$ is the integral 
$mk\times {\ell}k$ matrix $J^{\pi\circ \sigma}$. 
Noting that $\ker \epsilon\otimes \Z[G/K] =\Z^{k-1}$ finishes the proof. 
\end{proof}

To compute the abelianization of $K$ by Fox's method, one needs 
to row-reduce the matrix $J^{\pi\circ \sigma}$. In practical terms, 
this can be difficult, due to the rather big size of this matrix.  
If $K$ is a normal subgroup of $G$, a more efficient method 
is to first decompose the regular representation of $\G=G/K$ 
into irreducible representations.  Such a method will be 
described in Theorem~\ref{thm:genlibsak}.

\subsection{Representations of finite groups}
\label{subsec:reps}
Before proceeding, we need some basic facts from the representation 
theory of finite groups (see \cite{CR} as a reference).

\begin{definition} 
\label{def:sufflarge}
Let $\Gamma$ be a finite group, of order $\abs{\G}$.
A field $\K$ is {\em sufficiently large} with respect to $\G$
if the following two conditions are satisfied:

\begin{romenum}
\item \label{i}
The characteristic of $\K$ is $0$, or coprime to $\abs{\G}$.

\item \label{ii}
The field $\K$ contains all the $e$-roots
of unity, where $e$ is the exponent of $\Gamma$.

\end{romenum}
\end{definition}

Condition \eqref{ii} is satisfied if, for example, 
$\K$ is algebraically closed.  
If $\G=\Z_p$ is a cyclic group of prime order, 
and $q$ is a prime different from $p$, 
a sufficiently large field is $\K=\F_{q^s}$,
the Galois field of order $q^s$, 
where $s=\ord_p(q)$ is the least positive integer such that $p \mid (q^s-1)$.

If condition~\eqref{i} holds, then the group algebra $\K\Gamma$,
viewed as the regular representation of $\Gamma$, completely
decomposes into irreducible representations (Maschke).
If condition~\eqref{ii} holds, then $\K$ is a splitting field for 
$\G$ (Brauer).
Thus, if $\K$ sufficiently large, the regular
representation $\K\Gamma$ decomposes into (absolutely) irreducible
representations:
\begin{equation}
\label{eq:decomp}
\K\G=\bigoplus_{\rho\in Z}\bigoplus_{n_{\rho}} W_{\rho},
\end{equation}
where $Z=\Irrep(\G,\K)$ is the set of isomorphism classes of
irreducible $\K$ representations of $\Gamma$, and
$n_{\rho}$ is the dimension of the representation
$\rho:\G\to \GL(W_{\rho})$.  In particular,
$\abs{\Gamma}=\sum_{\rho\in Z} n_{\rho}^2$.

\subsection{Mod $q$  Betti numbers}
\label{subsec:modpbetti}
Let $b_1(G)=\rank H_1(G)=\dim_{\Q} H_1(G; \Q)$ be the 
first Betti number of $G$.  We shall write $b_1^{(0)}=b_1$.  
For a prime $q$, set 
\begin{equation*}
\label{def:thetaq}
b_1^{(q)}(G)=\dim_{\,\F_q}H_1(G; \F_q). 
\end{equation*}
Since homology commutes with direct sums, we have
$b_1^{(q)}(G)=\dim_{\,\K}H_1(G; \K)$, for any field
$\K$ of characteristic $q$.   We will call $b_1^{(q)}(G)$ the 
{\em mod $q$ (first) Betti number of $G$}.  

\begin{theorem}
\label{thm:genlibsak}
Let $G$ be a finitely-presented group.
Let $\lambda:G\surj\Gamma$ be a representation of $G$ onto
a finite group $\Gamma$.  Let $\K$ be a field, sufficiently
large with respect to $\G$.  Then,
if $K_{\l}=\ker \l$ and $q=\ch \K$:
\begin{equation}
\label{eq:genlibsak}
b_1^{(q)}(K_{\l}) = b_1^{(q)}(G)+\sum_{\rho\ne \triv}
n_{\rho}( \corank J^{\rho\circ \lambda}- n_{\rho}),
\end{equation}
where $J^{\rho\circ \lambda}$ is the Jacobian matrix of $G$,
followed by the representation $\rho\circ \lambda: G\to \GL(n_{\rho},\K)$,
and $\rho$ runs through all non-trivial, irreducible $\K$-representations
of $\Gamma$.
\end{theorem}

\begin{proof}
Let $G=\langle x_1,\dots ,x_{\ell}\mid r_1,\dots ,r_m\rangle$
be a finite presentation.
Since $\ch \K=q$, we have
$b_1^{(q)}(K_{\lambda})=\dim_{\K} H_1(K_{\lambda};\K)$.
As in the proof of Fox's Theorem \ref{thm:fox}, $H_1(K_{\l},\K)$ is
the first homology of the chain complex
\begin{equation}
\label{eq:chain}
(\K\G)^m \xrightarrow{J^{\l}} (\K\G)^{\ell}
\xrightarrow{\partial_1^{\l}} \K\G\xrightarrow{\epsilon}\K\to 0,
\end{equation}
obtained by tensoring \eqref{eq:chaincomplex} with $\K\G$
(viewed as a $G$-module via the representation $\l:G\to \G$).
In other words,
$b_1^{(q)}(K_{\l})=\dim_{\K} \ker\partial_1^{\l} -\dim_{\K} \im J^{\l}$.

Let $Z=\Irrep(\G,\K)$.  
In view of~\eqref{eq:decomp}, the chain complex~\eqref{eq:chain}
decomposes as:
\begin{equation*}
\label{eq:bigchain}
\bigoplus_{\rho\in Z}\bigoplus_{n_{\rho}} W_{\rho}^{m}
\xrightarrow{\oplus_{\rho}\oplus_{n_{\rho}} J^{\rho\circ\l}}
\bigoplus_{\rho\in Z}\bigoplus_{n_{\rho}} W_{\rho}^{\ell}
\xrightarrow{\oplus_{\rho}\oplus_{n_{\rho}} \partial_1^{\rho\circ\l}}
\bigoplus_{\rho\in Z}\bigoplus_{n_{\rho}} W_{\rho}
\xrightarrow{\epsilon}\K\to 0.
\end{equation*}
The twisted Jacobian matrices  $J^{\rho\circ\l}$ are of size
$m n_{\rho}\times \ell n_{\rho}$, and have entries in $\K$.
By Theorem~\ref{thm:fox},  $\rank J^{\triv\circ\l} =\ell - b_{1}^{(q)}(G)$.
Hence,
$\dim_{\K} \im J^{\l}=\ell-b_{1}^{(q)}(G)+\sum_{\rho\ne\triv}
n_{\rho}\rank J^{\rho\circ \lambda}$.  

Now notice that $\partial_1^{\l}$ has rank
equal to $\dim_{\K}\ker\epsilon=\abs{\G}-1$.  Hence,
$\dim_{\K} \ker \partial_1^{\l} ={\ell}\abs{\G}-(\abs{\G}-1)$.
Therefore,
\begin{equation*}
\label{eq:thetachain}
b_1^{(q)}(K_{\l})
={\ell}\abs{\G}-\abs{\G}+1-
\big(\ell-b_{1}^{(q)}(G)+
\sum_{\rho\ne\triv}n_{\rho} \rank J^{\rho\circ \lambda}  \big).
\end{equation*}
Since $\abs{\Gamma}=1+\sum_{\rho\ne\triv} n_{\rho}^2$,
we get formula \eqref{eq:genlibsak}.
\end{proof}

Notice that $\corank J^{\rho\circ\l}\ge 
\rank \partial_1^{\rho\circ\l}=n_{\rho}$.  
Hence, each term in the sum \eqref{eq:genlibsak} is non-negative, 
and so 
\begin{equation}
\label{eq:lowbound}
b_1^{(q)}(K_\l) \ge b_1^{(q)}(G).
\end{equation}
In view of this inequality, we are led to the following definition.
\begin{definition}
\label{def:betagammainv}
Let $G$ be a finitely-presented group, and let $\G$ be a finite group.
For $q=0$, or $q$ a prime not dividing $\abs{\G}$,
and $d$ a non-negative integer, put
\begin{equation*}
\label{eq:betagammadef}
\b_{\G,d}^{(q)}(G) :=  \#\big\lbrace K \lhd G \mid G/K\cong \G
\quad\text{and}\quad b_1^{(q)}(K) = b_1^{(q)}(G) + d \:\big\rbrace.
\end{equation*}
\end{definition}

In other words, $\b_{\G,d}^{(q)}(G)$ counts those 
normal subgroups of $G$, with factor group $\G$, 
for which the mod $q$ first Betti number
jumps by $d$, when compared to that of $G$.
Notice that:  
\begin{equation}
\label{eq:betagammasum}
\sum_{d\ge 0} \beta_{\G,d}^{(q)}(G)=\delta_{\G}(G).
\end{equation}

\subsection{Homology of finite abelian covers}
\label{subsec:abelian}

We may further refine Theorem~\ref{thm:genlibsak} in the case
when the group $\G$ is abelian.   We start with an immediate corollary.

If $\G$ is abelian, then all its irreducible representations over a
sufficiently large field $\K$ are $1$-dimensional.
Hence, taking $\K=\C$ in the above theorem, we obtain:

\begin{cor}[Libgober \cite{Li1}, Sakuma \cite{Sa}, Hironaka~\cite{Hr}]
\label{thm:libsakhir}
Let $\lambda:G\surj \Gamma$ be a representation of a
finitely-presented group $G$ onto a finite abelian
group $\Gamma$.  If $K_{\l}=\ker \l$, then:
\begin{equation}
\label{eq:libsak}
b_1(K_{\l}) = b_1(G) +\sum_{\rho\ne\triv}
( \corank J^{\rho\circ \lambda}- 1),
\end{equation}
where $\rho$ runs through all non-trivial, irreducible,
complex representations of $\Gamma$.
\end{cor}

For a representation $\rho:\G\to \K^*$,
let $\langle\rho\rangle$ be the (cyclic) subgroup of $\Hom(\G,\K^*)$
generated by $\rho$, and set  $m_{\rho}=\abs{\langle\rho\rangle}$.
Let $Z^{\wedge}=\Irrep^{\wedge}(\G,\K)$ be a set of representatives for the non-trivial,
irreducible $\K$-representations of $\G$, under the equivalence relation
$\rho_1\sim\rho_2\Longleftrightarrow 
\langle\rho_1\rangle=\langle\rho_2\rangle$.

\begin{theorem}
\label{thm:equivreps}
Let $\lambda:G\surj \Gamma$ be a representation of a
finitely-presented group $G$ onto a finite abelian
group $\Gamma$. Let $\K$ be a sufficiently
large field, of characteristic $q$.  Then:
\begin{equation}
\label{eq:grouped}
b_1^{(q)}(K_{\l}) = b_1^{(q)}(G)+\sum_{\rho\in Z^{\wedge}}
m_{\rho}( \corank J^{\rho\circ \lambda}-1).
\end{equation}
\end{theorem}

\begin{proof}
Assume $\rho_1\sim\rho_2$.
Let $C$ be the cyclic group generated by $\rho_1$. Then there is
an automorphism $\psi:C\to C$ such that $\psi(\rho_1)=\rho_2$. The
linear extension $\psi:\K C\to \K C$ is an isomorphism, taking
$J^{\rho_1\circ \l}$ to $J^{\rho_2\circ \l}$.  Consequently,
the $\K$-modules presented by these two matrices are isomorphic.
Hence,
$\corank J^{\rho_1\circ \lambda}=\corank J^{\rho_2\circ \lambda}$,
and so the contributions of $\rho_1$ and $\rho_2$ to the
sum \eqref{eq:genlibsak} are equal.
\end{proof}

The above theorem permits us to derive bounds and congruences 
on the mod~$q$ Betti numbers of normal subgroups  $K\nor G$ with 
$\G=G/K$ finite abelian, provided $q\nmid \abs{\G}$.

\begin{cor}
\label{cor:boundcong}
Let $G$ be a finitely-presented group, and $K$ a normal subgroup with 
$G/K$ finite abelian. Suppose $q=0$, or $q$ is a prime not dividing $k=\abs{G/K}$. 
\begin{enumerate}
\item \label{bound}
Let $\ell(G)$ be the minimal number of generators in a finite 
presentation for $G$.  Then:
\begin{equation*}
\label{eq:hbounds}
b_1^{(q)}(G) \le b_1^{(q)}(K) \le 
b_1^{(q)}(G)+(k-1)(\ell(G)-1).
\end{equation*}
\item \label{cong}
Let $p_1,\dots , p_r$ be the prime factors of $k$, and set
$D=\gcd(p_1-1,\dots,p_r-1)$. Then:
\begin{equation*}
\label{eq:cong}
b_1^{(q)}(K) \equiv b_1^{(q)}(G) \mod D.
\end{equation*}
\end{enumerate}

\end{cor}

\begin{proof}
\eqref{bound}  
The first inequality was already noted in \eqref{eq:lowbound}. 
To prove the second inequality, pick a finite presentation of 
$G$ with $\ell=\ell(G)$ generators, so that the 
Jacobian matrix $J=J_G$ has $\ell$ columns. 
For each representation $\rho\in Z^{\wedge}$,  
the matrix $J^{\rho\circ \lambda}$ 
also has $\ell$ columns, since $\rho$ is one-dimensional. 
Hence, $\corank J^{\rho\circ
\lambda}\le \ell$, for all $\rho\in Z^{\wedge}$.

\eqref{cong}
Each $m_{\rho}$ is of the form $\phi(r)$, for some
integer $r>1$ dividing $k$.  Hence, $D \mid m_{\rho}$,
for all $\rho\in Z^{\wedge}$.
\end{proof}

\section{Alexander matrices and characteristic varieties}
\label{sec:alexcv}

In this section, we introduce the characteristic varieties
of a finitely-presented group $G$, over an arbitrary field $\K$.  

\subsection{Alexander matrix}
\label{subsec:alexmat}
Let $G=\langle x_1,\dots ,x_{\ell}\mid r_1,\dots ,r_m\rangle$
be a finitely presented group.
Let $\phi:F_{\ell}\surj G$ be the presenting homomorphism, and
$\a:G \surj H_1(G)$ the abelianization map.
Fix an isomorphism $\chi: H_1(G)\xrightarrow{\cong} \Z^n 
\oplus\bigoplus_{i=1}^h
\Z_{e_i}^{n_i} $, where $e_i$ are distinct elementary divisors.
This identifies the group-ring $\Z H_1(G)$ with the ring
\begin{equation}
\label{eq:lambda}
\Lambda=\Z[t_1^{\pm 1},\dots,t_n^{\pm1}, s_{1,1},\dots ,s_{1,n_1},
\dots, s_{h,1},\dots, s_{h,n_h}]/
(s_{i,j}^{e_i}-1).
\end{equation}
The {\em Alexander matrix} of $G$ 
is the $\ell \times  m$ matrix with entries in $\Lambda$ given by
\begin{equation}
\label{eq:alexmat}
A_G=J_G^{\chi\circ \a},
\end{equation}
where recall 
$J_G=\begin{pmatrix}\frac{\partial r_i}{\partial x_j}\end{pmatrix}^{\phi}$
is the Fox Jacobian matrix associated to the given presentation of $G$.  
The {\em $d$-th Alexander ideal}, $E_d(A_G)$, is the ideal of
$\Lambda$ generated by the codimension~$d$ minors of $A_G$.
As is well-known, this ideal does not depend on the choice
of presentation for $G$ (but it does depend on the choice
of isomorphism $\chi$).

\subsection{Characteristic varieties}
\label{subsec:charvar}
Let $\K$ be a field,
and let $\K^*$ be its multiplicative group of units.
For $N$ a positive integer, let $\O_{N,\K}$ be the set of roots of unity
of order $N$ in $\K$.

Let $\Hom(G,\K^*)$ be the group of $\K$-valued characters of $G$.
The isomorphism $\chi:H_1(G)\to \Z^n \oplus\Z_{e_1}^{n_1} 
\oplus\cdots \oplus \Z_{e_h}^{n_h}$
identifies the character variety $\Hom(G,\K^*)$
with the product of affine algebraic tori
\begin{equation}
\label{eq:torus}
\T=(\K^*)^n \times (\O_{e_1,\K})^{n_1} \times \dots \times
(\O_{e_h,\K})^{n_h},
\end{equation}
viewed as a subset of the torus $(\K^*)^{n+n_1+\cdots +n_h}$.

\begin{definition}
\label{def:cv}
The {\em $d$-th characteristic variety} of the group $G$
(over the field $\K$) is the subvariety $V_d(G,\K)$ of the algebraic 
torus $\T=\Hom(G,\K^*)$, consisting of characters $\mathbf{t}:G\to\K^*$ such
that $f(\mathbf{t})=0$, for all $f \in E_d(A_G) \otimes \K$.
\end{definition}

In other words, $V_d(G,\K)$ is the $d$-th determinantal 
variety of $A_G\otimes \K$. As such, $V_d(G,\K)$ is also  
defined by the annihilator of the $d$-th exterior power 
of the Alexander module, $\coker A_G \otimes\K$, see \cite[pp.~511--513]{Ei}.  
If $d<\ell(G)$, this module has the same support as the Alexander invariant, 
$H_1(G',\K)$, where $G'$ is the commutator subgroup of $G$.  

The characteristic varieties of $G$ form a descending tower,
$\T=V_0\supseteq V_{1}\supseteq \cdots \supseteq V_{\ell(G)-1}\supseteq
V_{\ell(G)}=\emptyset$.  
This tower depends only on the isomorphism type of $G$,
up to a monomial change of basis in the algebraic torus
$(\K^*)^{n+n_1+\cdots +n_h}$.

\begin{remark}
\label{rem:hironaka}
The characteristic varieties $V_d(G,\K)$ may be
interpreted as the jumping loci for the cohomology 
of $G$ with coefficients in rank~$1$ local systems 
over $\K$.  More precisely, let 
\begin{equation}
\label{eq:jump}
\Sigma_d(G,\K)=\{\mathbf{t}\in \T \mid
\dim_{\K} H^1(G,\K_{\mathbf{t}}) \ge d\},
\end{equation}
where $\K_{\mathbf{t}}$ is the $G$-module $\K$ with
action given by the representation $\mathbf{t}:G\to \K^*$. 
Then, $V_d(G,\K)\setminus \{\mathbf{1}\}=
\Sigma_d(G,\K)\setminus \{\mathbf{1}\}$. 

For $\K=\C$, this was proved by Hironaka~\cite{Hr} 
(see also Libgober~\cite{Li2} and Cogolludo~\cite{Cog}).   
The proof given in \cite{Li2, Cog} can be adapted to work 
for an arbitrary field $\K$.  
Recall that $V_d(G,\K)$ is defined by  
$\ann\big(\bigwedge^d(H_1(G',\K))\big)$.  Thus, 
$\mathbf{t}\in V_d(G,\K) \Longleftrightarrow \dim_{\K}
H_1(G,\K_{\mathbf{t}})\ge d$ (here we need
$\mathbf{t}\ne \mathbf{1}$).   But $H_1(G,\K_{\mathbf{t}})\cong
H^1(G,\K_{\mathbf{t}^{-1}})$,  see \cite[p.~341]{Brd}, 
and we are done.
\end{remark}

\begin{remark}
\label{rem:resonance}
Closely related are the {\em resonance varieties} of the group $G$ 
(over the field $\K$), defined as 
$R_d(G,\K)=\{\l\in H^1(G,\K) \mid \dim_{\K} H^1(H^*(G,\K),\cdot \l) \ge d\}$,  
see \cite{Fa, LY, MS2}.  If all the relators of $G$ are commutators 
($r_j\in [F_{\ell},F_{\ell}]$, $\forall j$), 
then $R_d(G,\K)$ is the $d$-th determinantal variety of the 
{\em linearized} Alexander matrix of $G$, see \cite{MS2}.  
Moreover, as shown by Libgober~\cite{Li3}, the tangent cone at
$\mathbf{1}$  to $V_d(G,\C)$ is included in $R_d(G,\C)$. 
\end{remark}

\subsection{Depth of characters}
\label{subsec:depth}
Let $\mathbf{t}:G\to\K^*$ be a character. Since $\K^*$ is an abelian group,
$\mathbf{t}$ factors through the abelianization $\a:G \to H_1(G)$.
Let $A^{\mathbf{t}}:\K^m\to \K^{\ell}$
be the matrix obtained from $A$ by evaluating at $\mathbf{t}$.
(Under the isomorphism $\chi:\K H_1(G)\to\Lambda\otimes \K$,
the twisted Jacobian matrix $J^{\mathbf{t}}$ corresponds
to the twisted Alexander matrix $A^{\mathbf{t}}$.)
We then have:
\begin{equation}
\label{eq:twistalex}
  \mathbf{t}\in V_d(G,\K) \Longleftrightarrow \rank_{\K} A^{ 
\mathbf{t}}\le \ell-d-1.
\end{equation}

\begin{definition}
Let $G$ be a finitely-presented group, and $\K$ a field.
The {\em depth} of a character $\mathbf{t}:G\to \K^*$ (relative
to the stratification of $\Hom(G,\K^*)$ by the characteristic
varieties) is:
\[
d_{\K}(\mathbf{t})=\max\, \{ d \mid \mathbf{t}\in V_d(G,\K) \}.
\]
\end{definition}

Note that $0\le d_{\K}(\mathbf{t})\le \ell(G)-1$.  
Thus,  we can sharpen the upper bound from 
Corollary~\ref{cor:boundcong}\eqref{bound}, as
follows.  Let $K\nor G$ be a normal subgroup, 
with $G/K$ abelian of order $k$, and choose $\K$ to be
sufficiently large with respect to $G/K$.  Then: 
\begin{equation}
\label{eq:upperbound}
b_1^{(q)}(K) \le b_1^{(q)}(G)+(k-1) d_{\K}(G),
\end{equation}
where $d_{\K}(G)= \sup \{d_{\K}(\mathbf{t}) \mid
\mathbf{1}\ne \mathbf{t}\in \Hom(G,\K^*)\}$.

\section{Torsion points and Betti numbers}
\label{sec:bettitors}

Given a normal subgroup of $K\nor G$ with finite cyclic quotient, 
we interpret the first homology of $K$ with coefficients in a sufficiently 
large field $\K$, in terms of the stratification of the character 
variety $\Hom(G,\K^*)$ by the characteristic varieties of $G$.

\subsection{Homology of normal subgroups with cyclic quotient}
\label{subsec:cyclic}
Let $\Gamma=\Z_N$ be a finite cyclic group.
Let $\K$ be a field.  Assume that $\K$ is sufficiently
large with respect to $\Z_N$.  Then $\K$ contains
all the $N$-th roots of unity, and so there is a monomorphism
$\iota:\Z_N\inj\K^*$, sending a generator of $\Z_N$ to a
primitive $N$-th root of unity in $\K^*$.  Finally, for $j\ge 0$, let
$\psi_j:\K\to \K$ be the map $\psi_j(x)=x^j$.

\begin{theorem}
\label{thm:homcyclic}
Let $\l:G\surj \Z_N$ be a surjective homomorphism.
Let $\K$ be a field, sufficiently large with respect to $\Z_N$.
Set $K_{\l}=\ker(\l)$, and $q=\ch\K$.  Then
\begin{equation}
\label{eq:b1jump}
b_1^{(q)}(K_{\l}) = b_1^{(q)}(G)+\sum_{1\ne k\mid N} \phi(k) d_{\K}(\l^{N/k}),
\end{equation}
where $\l^{N/k}=\psi_{N/k}\circ \iota\circ \l$.
\end{theorem}

\begin{proof}
By Theorem~\ref{thm:equivreps}, we have
\begin{equation}
\label{eq:b1grouped}
b_1^{(q)}(K_{\l}) = b_1^{(q)}(G)+\sum_{\rho\in Z^{\wedge}}
m_{\rho}( \corank J^{\rho\circ \lambda}-1).
\end{equation}

Since $\K$ is sufficiently large, $\Hom(\Z_N,\K^*)\cong \Z_N$.
Let $\ord(\rho)$ be the order of a non-trivial 
representation $\rho:\Z_N\to \K^*$. 
It is readily seen that $\rho\sim \rho' \Longleftrightarrow \ord(\rho)=\ord(\rho')$. 
Thus, the assignment $\rho\mapsto \ord(\rho)$ establishes a bijection 
between $Z^{\wedge}=\Irrep^{\wedge}(\Z_N,\K)$ and the set 
of non-unit divisors of $N$. Moreover, $m_{\rho}=\phi(k)$, 
where $k=\ord(\rho)$.

Now consider the representation $\rho=\psi_{N/k}\circ \iota$.
Clearly, the order of $\rho$ is $k$. 
Let $\l^{\vee}:\Hom(\Z_N,\K^*)\to\Hom(G,\K^*)$ be the dual homomorphism.
We then have $\l^{\vee}(\rho)=\l^{N/k}$. Furthermore, by \eqref{eq:twistalex},
we have
\begin{equation}
\label{eq:lvee}
\l^{\vee}(\rho) \in V_d(G,\K) \Longleftrightarrow
\corank_{\K} A^{\rho\circ\l} \ge d+1.
\end{equation}
The conclusion follows at once.
\end{proof}

\begin{cor}
\label{cor:homprime}
Let $K\nor G$ be a normal subgroup of
prime index~$p$.  Write $K=\ker(\l:G\to \Z_p)$.
Let $q=0$, or $q$ a prime, $q\ne p$.
Let $\K$ be a field of characteristic $q$
which contains all the $p$-roots of unity---for example,
$\K=\C$, or $\K=\F_{q^s}$, where $s=\ord_p(q)$.
Then:
\begin{equation}
\label{eq:thetajump}
b_1^{(q)}(K) = b_1^{(q)}(G)+(p-1)d_{\K}(\l).
\end{equation}
\end{cor}

\subsection{Distribution of mod $q$ Betti numbers}
\label{subsec:dist}
In view of the above corollary, it makes sense to define
\begin{equation}
\label{eq:betadef}
\b_{p,d}^{(q)}(G) := \b_{\Z_p,(p-1)d}^{(q)}(G)
\end{equation}
In other words, $\b_{p,d}^{(q)}(G)$ counts those index~$p$,
normal subgroups of $G$ for which the mod~$q$ first Betti number
jumps by $(p-1)d$, when compared to that of $G$.
By \eqref{eq:betagammasum}, we have 
$\sum_{d\ge 0} \beta_{p,d}^{(q)}(G)=\delta_{\Z_p}(G)$.  
Hence, by Theorem~\ref{thm:abelcount}:
\begin{equation}
\label{eq:betasum}
\sum_{d\ge 0} \beta_{p,d}^{(q)}(G)=\frac{p^n-1}{p-1}, \quad 
\text{where $n=b_1^{(p)}(G)$}.
\end{equation}
For simplicity, we shall write sometimes 
$\beta_p^{(q)}=\big(\beta_{p,1}^{(q)},\dots,\beta_{p,k}^{(q)}\big)$,
if $\beta_{p,d}^{(q)}=0$, for $d>k$. Also, we will abbreviate
$\b_{p,d}=\b_{p,d}^{(0)}$.  Note that $\beta_{p,0}^{(q)}$
is determined from \eqref{eq:betasum} by the sequence $\beta_p^{(q)}$,
and the mod $p$ first Betti number of $G$.

Let 
\begin{equation}
\label{eq:torspd}
\Tors_{p,d}(G,\K)=\{  \mathbf{t}\in \Hom(G,\K^*) \mid
\mathbf{t}^p=\mathbf{1}\ \text{and}\ \mathbf{t}\ne \mathbf{1}\}\cap 
V_d(G,\K)
\end{equation}
be the set of characters on $V_d(G,\K)$ of order exactly equal to~$p$. 
As a direct consequence of Corollary~\ref{cor:homprime}, we obtain:
\begin{theorem}
\label{thm:torscount}
Let $G$ be a finitely-presented group, $p$ a prime, 
and $q=0$, or $q$ a prime, distinct from $p$.  
If $\K$ is a field of characteristic $q$ containing
all $p$-roots of unity, then:
\[
\b_{p,d}^{(q)}(G)=\frac{\abs{\Tors_{p,d}(G,\K)
\setminus \Tors_{p,d+1}(G,\K)}}{p-1}.
\]
In particular,
$\beta_{p,d}(G)=\tfrac{1}{p-1}\abs{\Tors_{p,d}(G,\C)
\setminus \Tors_{p,d+1}(G,\C)}$. Also, if $q>0$, then 
$\b_{p,d}^{(q)}(G)=\tfrac{1}{p-1}\abs{\Tors_{p,d}(G,\F_{q^s})
\setminus \Tors_{p,d+1}(G,\F_{q^s})}$, where $s=\ord_p(q)$. 
\end{theorem}

\begin{remark}
\label{rem:nuinv}
This result does not say anything about the 
distribution of mod~$p$ Betti numbers of index $p$ 
subgroups. Even so, there is a class of groups for which 
an analogous formula holds for $q=p$, with the characteristic varieties 
replaced by the resonance varieties (over the field $\F_p$). 
Indeed, let $G=\langle x_1,\dots , x_n \mid r_1,\dots ,r_m\rangle$ 
be a commutator-relators group, with $H_2(G)$ torsion-free, 
and let $Q=G/[G,[G,G]]$ be its 
second nilpotent quotient.  Set 
$\nu_{p,d}(Q)=\#\big\lbrace K \lhd Q \mid \ind{Q:K}=p
\ \text{and}\ b_1^{(p)}(K) = n + d \:\big\rbrace$.  Then, 
according to \cite[Theorem~4.19]{MS2}:
\begin{equation}
\label{eq:nupd}
\nu_{p,d}(Q)=\tfrac{1}{p-1}\abs{R_d(Q,\F_p)\setminus R_{d+1}(Q,\F_p)}.
\end{equation}
\end{remark}

\subsection{Computations of $\beta$-invariants}  
\label{subsec:exbeta}
We conclude this section with some sample computations 
of the invariants $\beta_{p,d}^{(q)}(G)$, for some familiar 
finitely-presented groups $G$. 

\begin{example}
\label{ex:freegp}
Let $G=F_n$ be the free group of rank $n$.  Evidently,  
$V_0(G,\K)=\cdots =V_{n-1}(G,\K)=(\K^{*})^n$, and $V_n(G,\K)=\{ \mathbf{1} \}$, 
for all $\K$. Hence, for all $q$:
\begin{equation*}
\label{eq:betafree}
\beta_{p,n-1}^{(q)}(G)=\frac{p^n-1}{p-1},\quad \text{and}\quad
\beta_{p,d}^{(q)}(G)=0, \text{ for $d\ne n-1$.}
\end{equation*}
\end{example}

\begin{example}
\label{ex:prodfree}
Let $G=F_m\times F_n$ ($m\ge n$) be the product of two free groups.  Then:
\begin{equation*}
V_d(G,\K)=
\begin{cases}
(\K^*)^{m+n}\quad &\text{if }\, d=0,\\
(\K^*)^{n}\cup (\K^*)^{m} \quad
&\text{if }\, 0<d<n\\
(\K^*)^{m} \quad &\text{if }\, n\le d < m,\\
\{ \mathbf{1} \}\quad &\text{if }\, d=m.
\end{cases}
\end{equation*}
where $(\K^*)^{n}=\{t_1=\cdots =t_m=1\}$ and $(\K^*)^{m}=\{t_{m+1}=\cdots
=t_{m+n}=1\}$ (see \cite{CScv}). Hence:
$\beta_{p,0}^{(q)}(G)=\tfrac{(p^m-1)(p^n-1)}{p-1}$,
$\beta_{p,n-1}^{(q)}(G)=\tfrac{p^n-1}{p-1}$,
$\beta_{p,m-1}^{(q)}(G)=\tfrac{p^m-1}{p-1}$,
and $\beta_{p,d}^{(q)}(G)=0$, otherwise.

More generally, consider the product $G=F_{n_1}\times \cdots \times F_{n_k}$, 
with $n_1\ge \cdots \ge n_k$. Write $\mathbf{n}=(n_1,\dots,n_k)$, and recall that
$\abs{\mathbf{n}}=\sum_{i=1}^{k} n_i$  and $m_d(\mathbf{n})=\#\{j\mid n_j=d\}$. 
We then have:
\[
\beta_{p,0}^{(q)}(G)=\frac{p^{\abs{\mathbf{n}}}-1-\sum_{i=1}^{k} 
(p^{n_i}-1)}{p-1}, 
\qquad 
\beta_{p,d-1}^{(q)}(G)=m_{d}(\mathbf{n})\frac{p^d-1}{p-1},\ \text{ 
for } d>1.
\]
\end{example}

\begin{example}
\label{ex:surfacegp}
Let $G=\pi_1(\#_{g} T^2)$
be the fundamental group of a closed, orientable surface of genus $g\ge 1$, 
with presentation
$G=\langle x_1,\dots ,x_{2g} \mid [x_1,x_2]\cdots [x_{2g-1},x_{2g}]=1\rangle$.
Then: $V_d(G,\K)=(\K^{*})^{2g}$, for $d<2g-1$, and
$V_{2g-1}(G,\K)=\{ \mathbf{1} \}$ (see \cite{Hr}).  Hence: 
\[
\beta_{p,2g-2}^{(q)}(G)=\frac{p^{2g}-1}{p-1}, \text{ and }
\beta_{p,d}^{(q)}(G)=0, \text{ for $d\ne 2g-2$.}
\]
\end{example}

\begin{example}
\label{ex:nonorient}
Let $G=\pi_1(\#_{n} \RP^2)$
be the fundamental group of a closed, non-orientable surface of genus $n\ge 1$, 
with presentation
 $G=\langle x_1,\dots , x_{n} \mid x_1^2\cdots x_n^2=1 \rangle$.
The isomorphism $\chi:H_1(G)\to \Z^{n-1}\oplus \Z_2$, given by $\chi(x_i)=t_i$
for $i<n$ and $\chi(x_1\cdots x_n)=s$, identifies $\Z H_1(G)$ with 
$\Z[t_1^{\pm 1},\dots,t_{n-1}^{\pm1}, s]/(s^2-1)$.
The Alexander matrix $A_G$ is
\[
\begin{pmatrix}
1+t_1 & t_1^2(1+t_2) & \cdots &t_1^2\cdots t_{n-2}^2 (1+t_{n-1}) &
t_1^2 \cdots t_{n-1}^2(1+t_1^{-1}\cdots t_{n-1}^{-1}s)
\end{pmatrix}.
\]
If $\ch \K\ne 2$, the character variety $\Hom(G,\K^*)$ is isomorphic
to $\T=(\K^*)^{n-1}\times \{\pm 1\}$, and the characteristic varieties are:
$V_0=\dots =V_{n-2}=\T$, and $V_{n-1}=\{(-1,\dots,-1,(-1)^n)\}$.
If $\ch \K=2$, then $V_0=\dots =V_{n-2}=(\K^*)^{n-1}$, and
$V_{n-1}=\{\mathbf{1}\}$.  Hence:
\[
\beta_{2,n-2}^{(q)}(G)=2^n-2, \quad
\beta_{2,n-1}^{(q)}(G)=1, \quad
\beta_{p,n-2}^{(q)}(G)=\frac{p^{n-1}-1}{p-1},
\]
and $\beta_{p,d}^{(q)}(G)=0$, otherwise.
\end{example}

\section{Counting metabelian representations}
\label{sec:metacyclic}

We now return to the Hall invariants of a finitely-presented group $G$.
We show how to compute $\delta_{\G}(G)$, for the split metabelian groups
$\G=\Z_q^s\rtimes \Z_p$, in terms of torsion points on the
characteristic varieties of $G$, over the Galois field $\F_{q^s}$.

\subsection{A class of metabelian groups}
\label{subsec:metagp}
For two distinct primes $p$ and $q$, we define the
metabelian group $M_{p,q^s}$ to be the (non-trivial) split extension
\begin{equation}
M_{p,q^s}=\Z_q^s \rtimes_{\s} \Z_p=\langle a_1, \dots,a_s, b\mid
a_i^q =[a_i,a_j]=b^p=1, b^{-1}a_i
b=\s(a_i)\rangle,
\end{equation}
where $s=\ord_p(q)$ is the order of $q\bmod p$ in $\Z_p^*$,
and $\sigma$ is an automorphism of $\Z_q^s$, of order exactly $p$.

Note that $\s$ must act trivially on any proper,
invariant subgroup of $\Z_q^s$, and so, all proper
subgroups of $M_{p,q^s}$ are abelian.  Implicit in the
definition is the assertion that such automorphism $\s$ exists,
and that the isomorphism type of $M_{p,q^s}$ does not depend
on its choice.  This is proved in the following lemma,
which also gives the order of the automorphism group
of $M_{p,q^s}$.

\begin{lemma}
\label{lem:mpq}
Let $p$ and $q$ be distinct primes, and let $s=\ord_p(q)$.  Then:
\begin{enumerate}
\item \label{exist}
There exists an automorphism  $\s\in\Aut(\Z_q^s)$ of order $p$.
\item \label{unique}
If $\psi$ is another automorphism of $\Z_q^s$ of order $p$,
then $\Z_q^s\rtimes_{\s} \Z_p\cong \Z_q^s\rtimes_{\psi} \Z_p$.
\item \label{auto}
$\abs{\Aut(\Z_q^s\rtimes_{\s} \Z_p)}=sq^s(q^s-1)$.
\end{enumerate}
\end{lemma}

\begin{proof}
\eqref{exist}
The cyclotomic polynomial $Q_p=t^{p-1}+\cdots+t+1$ factors over
the field $\F_q$ into $(p-1)/s$ distinct, monic,
irreducible polynomials of degree $s$.  If $f$ is any one of those
factors, then the field $\F_{q^s}$ is isomorphic to $\F_q[t]/(f)$.
Let $\s$ be the automorphism of $\F_q[t]/(f)$
induced by multiplication by $t$ in $\F_q[t]$.
Clearly, $\s$ has order $p$.

Note the following: If we view $\Z_q^s$ as the $\F_q$-vector space
with basis $\{1,t,\dots,t^{s-1}\}$, then
$\s\in\Aut(\Z_q^s)\cong\GL(s,q)$ may be identified
with the companion matrix of $f$, and so $f=f_{\s}$,
the characteristic polynomial of $\s$.
Alternatively, if we view $\Z_q^s$ as the additive group of the
field $\F_{q^s}=\F_q(\xi)$, where $\xi$ is a primitive
$p$-th root of unity, then $\s=\cdot \xi\in\Aut(\F_q(\xi))$.

\eqref{unique}
Notice that $\Z_p^s\rtimes_{\s}\Z_p$ is isomorphic to
$\Z_p^s\rtimes_{\s^l}\Z_p$, for any $0<l<p$:
The mapping $a_i \mapsto a_i$, $b \mapsto b^k$,
where $k=l^{-1}$ in the multiplicative group $\Z_p^*$,
provides such an isomorphism.

Now let $\psi$ be an arbitrary matrix of order $p$ in $\GL(s,q)$.
The characteristic polynomial of $\psi$ must be one of the $(p-1)/s$
irreducible factors of $Q_{p}$.  All such factors are the characteristic
polynomials of some power of $\s$.  Thus, $f_{\psi}=f_{\s^{l}}$,
for some $0<l<p$. An exercise in linear algebra shows that
there is a matrix $\phi\in \GL(s,q)$ such that $\psi=\phi\s^l\phi^{-1}$.
Hence, $\Z_q^s\rtimes_{\psi} \Z_p\cong \Z_q^s\rtimes_{\s^l} \Z_p$.

\eqref{auto}
Let $\Phi\in\Aut(\Z_q^s \rtimes_{\s} \Z_p)$.  Every element in the semi-direct product
$\Z_q^s \rtimes_{\s} \Z_p$ has unique normal form $u b^k$, for some 
$u\in \Z_q^s$ and $0\le k \le p-1$.  Write $\Phi(b)=v b^{l}$.   
Straightforward computations show that $\Phi$ leaves 
the subgroup $\Z_q^s$ invariant, and that the restriction 
$\phi=\Phi\lvert_{\Z_q^s}$ satisfies $\phi\s\phi^{-1}=\s^l$.   
The number of solutions $\phi\in\GL(s,q)$ of this equation is $q^s-1$. 

The count of automorphisms of $M_{p,q^s}$ then follows from the 
following claim:
There are precisely $s$ values $1\le l \le p-1$ for which
$f_{\s^l}=f_{\s}$.
To prove the claim, notice that $\prod_{l=1}^{p-1} f_{\s^l} =Q_p^s$
(both polynomials factor completely over $\F_{q^s}$
into linear factors, and the factorizations coincide).
But, over $\F_q$, the polynomial $Q_p$ has  $(p-1)/s$ distinct
irreducible factors, all with multiplicity one, and so $Q_p^s$ has $(p-1)/s$
irreducible factors, each appearing exactly $s$ times.
\end{proof}

\begin{example}
If $p| (q-1)$, then $s=1$, and $\s:\F_q\to \F_q$ is given by
$\s(1)=r$, where $r^p=1\ (\bmod\, q)$ and $r\ne 1$. Thus,
$M_{p,q}=\langle a,b\mid a^p =b^q=1, a^{-1}ba=b^r\rangle$
is the  metacyclic group of order $pq$, and $\Aut(M_{p,q})\cong M_{q-1,q}$.
Well-known examples are the dihedral groups $D_{2q}=M_{2,q}$,
and, in particular, the symmetric group $S_3=D_6$.

If $p=3$ and $q=2$, then $s=2$ and
$\s=\bigl(\begin{smallmatrix}0&1\\1&1\end{smallmatrix}\bigr)$.
The group $M_{3,4}=\Z_2^2\rtimes_{\s} \Z_3$ is isomorphic to
the alternating group $\Alt_4$, and $\Aut(M_{3,4})\cong S_4$.
\end{example}

\subsection{Metabelian representations}
\label{subsec:metareps}
We now study the homomorphisms from a finitely-pre\-sented group $G$
to the metabelian group $\Gamma=M_{p,q^s}$.  Our approach is modelled
on that of Fox~\cite{Fx2}, where similar results are obtained
in the case where $G$ is a link group (with the Wirtinger
presentation), and $\Gamma=M_{p,q}$ is metacyclic.

Let $\phi:F_{\ell}\surj G$ be a presenting homomorphism,
with $F_{\ell}=\langle x_1,\dots, x_{\ell}\rangle$.
Let $\G=B\rtimes_{\s} C$ be a semidirect product of abelian groups,
with monodromy homomorphism $\s:C\to \Aut(B)$. Denote also by $\s$
the linear extension to group-rings, $\s: \Z C\to \End(B)$.
Finally, let $\rho:G\to C$ be a homomorphism, and set
$\bar{\rho}=\rho\circ\phi: F_{\ell}\to C$.

\begin{lemma}
\label{lem:abword}
Suppose $\l:F_{\ell}\to \G$ is a lift of $\bar\rho$,
given on generators by $\l(x_i)=b_i\bar{\phi}(x_i)$.
If $\l(w)=b\bar{\rho}(w)$, then the following equality holds in $B$:
\[
b=\sum_{j=1}^{\ell} \s\bar{\rho}\big(\tfrac{\partial w}{\partial 
x_j}\big)(b_j).
\]
\end{lemma}

\begin{proof}
The proof is by induction on the length of the word $w\in F_{\ell}$.
If $w=1$, the equality holds trivially.
Suppose $w=u x_i^{e}$, where $e=\pm 1$.  Put $\l(u)=b'\bar{\rho}(u)$.
We then have: $b\bar{\rho}(w)=\l(w)=\l(u 
x_i^{e})=b'\bar{\rho}(u)(b_i\bar{\rho}(x_i))^e$.
Rewriting this last word in normal form (in the semidirect product
$\G=B \rtimes_{\s} C$), we obtain the following equality
(in the additive group $B$):
\begin{eqnarray}
\label{eq:beta}
b=b'+e\, \s\bar{\rho}\big(u x_i^{(e-1)/2}\big)(b_i).
\end{eqnarray}

Taking Fox derivatives of $w=u x_i^{e}$, and applying $\bar{\rho}$, gives:
\begin{equation}
\bar{\rho}\big(\tfrac{\partial w}{\partial x_j}\big)=
\bar{\rho}\big(\tfrac{\partial u}{\partial x_j}\big) +
e\, \bar{\rho}\big(u x_i^{(e-1)/2}\big)\delta_{ij}.
\end{equation}
Now apply $\s$, evaluate at $b_j$, and sum over $j$:
\begin{align*}
\sum_{j=1}^{\ell} \s\big(\bar\rho\big(
\tfrac{\partial w}{\partial x_j} \big)\big) (b_j)
&=
\sum_{j=1}^{\ell} \s\big(\bar\rho\big(
\tfrac{\partial u}{\partial x_j}\big)\big) (b_j)+
	e\, \s\big(\bar{\rho}\big(u x_i^{(e-1)/2}\big)\big) (b_i) &&  \\
     &= b'+ e\, \s\big(\bar{\rho}\big(u x_i^{(e-1)/2}\big)\big) (b_i)
\quad \text{by induction hypothesis}&\\
     &= b \quad\text{by ~\eqref{eq:beta}.}
&\quad\qed
\end{align*}
\renewcommand{\qed}{}
\end{proof}

For $p$ and $q$ distinct primes, with $s=\ord_p(q)$,
let $M_{p,q^s}=\Z_{q}^s\rtimes_{\s} \Z_p$
be the split metabelian group defined in 
\ref{subsec:metagp}. 
Let $b$ be a generator of the cyclic group $\Z_p$.
Viewing $\Z_q^s$ as the additive group of the field $\K=\F_q(\xi)$,
where $\xi\in \K^*$ is a primitive $p^{\text{th}}$ root of unity,
we may take $\s(b)=\cdot \xi\in\Aut(\K)$.  In particular,
this identifies $\Z_p$ as a subgroup of $\K^*$, and thus,
$\Hom(G,\Z_p)$ as a subset of the character variety $\Hom(G,\K^*)$.

\begin{prop}
\label{prop:metabelcount}
The number of homomorphisms, respectively epimorphisms from the
finitely-presented group $G$ to the metabelian group $M_{p,q^s}$
is given by
\begin{align*}
\abs{\Hom(G,M_{p,q^s})}&=\sum_{\rho\in\Hom(G,\Z_p)} q^{sd_{\K}(\rho)+s},
\\
\abs{\Epi(G,M_{p,q^s})}&=\sum_{\mathbf{1}\ne\rho\in\Hom(G,\Z_p)} q^s 
(q^{sd_{\K}(\rho)}-1) .
\end{align*}
where $\K=\F_{q^s}$, and the sums are over
representations $\rho:G\to \Z_{p}\subset \K^*$.
\end{prop}

\begin{proof}
Let $G=\langle x_1,\dots ,x_{\ell}\mid r_1,\dots, r_{m}\rangle$ be a
presentation for $G$.  Let $\rho:G\to \Z_p$ be a representation,
given by $\rho(x_i)=b^{\b_i}$.  We want to lift it to a
representation $\l:G\to M_{p,q^s}$.  Such a representation
is given by $\l(x_i)=u_ib^{\b_i}$,
where $u_i\in \Z_q^s$.  In view of Lemma~\ref{lem:abword},
we must solve the following system of
equations over $\K=\F_q(\xi)$:
\begin{equation}
\label{eq:metafoxstar}
\sum_{j=1}^{\ell} A_{k,j}(\xi^{\b_1}, \dots, \xi^{\b_n})\cdot u_j
= 0,\quad 1\le k\le m,
\end{equation}
where $A_G=(A_{k,j})$ is the Alexander matrix of $G$.
This system has $q^{sd_{\K}(\rho)+s}$ solutions.
Starting now with a non-trivial representation $\rho:G\to \Z_p$,
all such solutions give rise to surjective representations $\l:G\to 
M_{p,q^s}$,
except $q^s$ of them, which give rise to abelian representations.
\end{proof}

This proposition, together with Lemma~\ref{lem:mpq}\eqref{auto},
imply the following.

\begin{theorem}
\label{thm:delmpqs}
If $p$ and $q$ are distinct primes, and $s=\ord_p(q)$, then:
\begin{equation*}
\label{eq:dmpq}
\delta_{M_{p,q^s}}(G) =\frac{p-1}{s(q^s-1)} \sum_{d\ge 1}
\b_{p,d}^{(q)}(G) (q^{sd}-1).
\end{equation*}
\end{theorem}

\begin{example}
\label{ex:delmpqfree}
For free groups, Theorem~\ref{thm:delmpqs} gives:
\begin{equation}
\label{eq:hallfree}
\delta_{M_{p,q^s}}(F_n) = \frac{(p^{n}-1)(q^{s(n-1)}-1)}{s(q^s-1)},
\end{equation}
since $\b_{p,n-1}^{(q)}(F_n)=\frac{p^n-1}{p-1}$,
and the other terms in the sum vanish.
In particular, this recovers formulas \eqref{eq:deltameta}
(when $p\mid q-1$) and \eqref{eq:deltalt} (when $p=3$, $q=2$).
For a product of free groups, we get:
\begin{equation}
\label{eq:hallprod}
\delta_{M_{p,q^s}}(F_{n_1}\times \cdots \times F_{n_k})= \frac{\sum_{i=1}^{k}
(p^{n_i}-1) (q^{s  (n_i-1))}-1)}{s(q^s-1)}.
\end{equation}
\end{example}

\begin{example}
\label{ex:delmpqsurf}
For orientable surface groups of genus $g\ge 1$, Theorem~\ref{thm:delmpqs} gives:
\begin{equation}
\label{eq:hallsurf}
\delta_{M_{p,q^s}}(G)= \frac{(p^{2g}-1)(q^{2s(g-1)}-1)}{s(q^s-1)}.
\end{equation}
For non-orientable surface groups of genus $n\ge 1$, we get:
\begin{align}
\label{eq:hallnonorient}
\delta_{D_{2q}}(G)&=\frac{(2^{n}-2)(q^{n-2}-1)+q^{n-1}-1}{q-1},\\
\delta_{M_{p,q^s}}(G)& = \frac{(p^{n-1}-1)(q^{s(n-2)}-1)}{s(q^s-1)}.
\end{align}
\end{example}

Table \ref{tab:hall2} gives the values of $\delta_{\Gamma}(G)$, 
for some of the groups $G$ in Examples~\ref{ex:delmpqfree} and \ref{ex:delmpqsurf}, 
and for some finite groups $\G$ of small order. 

\renewcommand{\arraystretch}{1.15}
\begin{table}
\[
\footnotesize
\begin{array}
{|l||r|r|r|r|r|r|r|r|r|}
\hline
\qquad G\ \backslash\ \Gamma
& \Z_2 & \Z_3 & \Z_2^2 & \Z_4  & \Z_2\oplus\Z_4 &  \Z_8 & S_3
& \Alt_4 & M_{3,7}\\
\hline\hline
F_2   & 3 &  4 & 1  &  6 &3 & 12 &  3 & 4  &    8 \\
F_3   & 7 & 13 & 7  & 28   &42 &  112& 28 & 65  &  208\\
F_4   & 15& 40 & 35 &120  &420 &  960&195& 840  &4,560\\
\hline
F_2\times F_1    & 7 & 13  & 7  & 28 &    42 &   112 &  3 &  4 &   8 \\
F_2\times F_2    & 15& 40  & 35 &120 &   420 &   960 &  6 &  8 &  16 \\
F_3\times F_1    & 15& 40  & 35 &120 &   420 &   960 & 28 & 65 & 208 \\
F_3\times F_2    & 31& 121 &155 &496 & 3,720 & 7,936 & 31 & 69 & 216 \\
\hline
\pi_1(\#_2 T^2)  &  15 &    40 &     35 &    120 &    420    & 960      
&     60  &      200 &         640\\
\pi_1(\#_3 T^2)      &  63 &   364 &    651 &  2,016 &  31,248   & 64,512    
& 2,520  &    30,940 &     291,200\\
\pi_1(\#_4 T^2)  & 255 & 3,280 & 10,795 & 32,640 & 2,072,640 & 4,177,920 
& 92,820 & 4,477,200 & 128,628,480\\
\hline
\pi_1(\#_2 \RP^2) & 3&  1 &   1 &   2 &     1&    2 &  1 &  0 &     0 \\
\pi_1(\#_3 \RP^2) & 7&  4 &   7 &  12 &    18&    24& 10 &  4 &     8\\
\pi_1(\#_4 \RP^2) &15& 13 &  35 &  56 &   196&   224& 69 & 65 &   208\\
\pi_1(\#_5 \RP^2) &31& 40 & 155 & 240 & 1,800& 1,920& 430&840 & 4,560\\
\hline
\end{array}
\]
\\[6pt]
\caption{\textsf{$\Gamma$-Hall invariants of some 
finitely presented groups $G$.}}
\label{tab:hall2}
\end{table}
\renewcommand{\arraystretch}{1.1}

\section{Counting finite-index subgroups}
\label{sec:lis}
We now discuss some other invariants of a
finitely-generated group $G$, obtained by counting
finite-index subgroups of $G$ in various ways.  
If $G$ is finitely-presented, and the index is 
low, these invariants can be computed 
from the characteristic varieties of $G$, and 
some simple homological data.

\subsection{Subgroups of finite index}
\label{subsec:findex}
For each positive integer $k$, let
\begin{equation}
a_{k}(G)=\text{number of index $k$ subgroups of $G$}.
\end{equation}
Also, let $h_l(G)=\sigma_{S_l}(G)$ be the number of homomorphisms 
from $G$ to the symmetric group $S_l$.  The following well-known formula of
Marshall~Hall~\cite{Ha} (see also \cite{Lu}) computes $a_k$ in terms of 
$h_1, \dots , h_k$ (starting from  $a_1=h_1=1$):
\begin{equation}
\label{eq:hall}
a_{k}(G)=\frac{1}{(k-1)!}\, h_{k}(G) -
\sum_{l=1}^{k-1} \frac{1}{(k-l)!}\, h_{k-l}(G) a_{l}(G).
\end{equation}

For the free group $G=F_n$, we have $h_k(F_n)=(k!)^n$, and so, as noted 
by M.~Hall, 
\begin{equation}
\label{eq:akfree}
a_k(F_n)=k(k!)^{n-1}-\sum_{l=1}^{k-1}((k-l)!)^{n-1} a_l(F_n).
\end{equation}

For the free abelian group $G=\Z^n$, a result of Bushnell and Reiner~\cite{BuR} 
gives $a_k(\Z^n)$ recursively, starting from $a_k(\Z)=1$:
\begin{equation}
\label{eq:bushreiner}
a_k(\Z^n)=\tsum_{d|k} a_d(\Z^{n-1})\left(\frac{k}{d}\right)^{n-1}
\end{equation}
(see \cite{Lin} for a simple proof, using the Hermite normal 
form of integral matrices).   Equivalently, 
$\zeta_{\Z^n}(s) = \prod_{i=0}^{n-1} \zeta(s-i)$, 
where $\zeta_{G}(s)=\sum_{k=1}^{\infty} a_k(G) k^{-s}$ 
is the zeta function of the group $G$, and $\zeta(s)$ is 
the classical Riemann zeta function (see \cite{Lu} for 
a detailed discussion). 

For surface groups $G$, the numbers $a_k(G)$ were computed by 
Mednykh~\cite{Me}.

As an application of our methods, we express the number of index
$2$ and $3$ subgroups of a finitely-presented group, in terms of its
characteristic varieties.

\begin{theorem}
\label{thm:index2and3}
Let $G$ be a finitely-presented group. Set $n_p=b_1^{(p)}(G)$. 
Then, the number of index $2$ and $3$ subgroups of $G$ is given by:
\begin{align*}
a_2(G)&=2^{n_2}-1, \\
a_3(G)&=\tfrac{1}{2}(3^{n_3}-1) + \tfrac{3}{2} 
\sum_{d\ge 1} \beta_{2,d}^{(3)}(G) (3^d-1). 
\end{align*}
\end{theorem}

\begin{proof}
Clearly, $a_2=h_2-1=\delta_{\Z_2}$, and the first identity 
follows from Theorem~\ref{thm:abelcount}.

M.~Hall's formula \eqref{eq:hall} gives $a_3=\frac{1}{2}h_3-\frac{3}{2}h_2+1$. 
Recall that the subgroup lattice of the symmetric group $S_3=M_{2,3}$ is
$L(S_3)=\{1,\Z_2,\Z_2,\Z_2,\Z_3,S_3\}$ (see~Figure~\ref{fig:sgplattice}).  
Thus, P.~Hall's formula \eqref{eq:hallenum1} gives 
$h_3=1+3\delta_{\Z_2}+2\delta_{\Z_3}+6\delta_{S_3}$. 
Using  Theorem~\ref{thm:abelcount}, we get:
\begin{equation}
\label{eq:a3dels3}
a_3(G)=\tfrac{1}{2}(3^{n_3}-1) + 3\delta_{S_3}(G).
\end{equation}
The second identity now follows from Theorem~\ref{thm:delmpqs}.
\end{proof}

For example, $a_3(F_{n})=3(3^{n-1}-1)2^{n-1}+1$, which agrees
with M.~Hall's computation.  Also, 
$a_3(\pi_1(\#_{g}T^2))=3(3^{2g - 2} - 1)(2^{2g - 1} + 1) + 4$ 
and 
$a_3(\pi_1(\#_{n}\RP^2))=3(3^{n - 2} - 1)(2^{n - 1} + 1) + 4$, 
which agrees with Mednykh's computation.

\begin{remark}
To compute $a_4(G)$ by the same method, one needs to do more work. Indeed, 
$a_4=\frac{1}{6}h_4-\frac{2}{3}h_3-\frac{1}{2}h_2^2+2h_2-1$, and 
\begin{equation}
h_4=1+9\delta_{\Z_2}+8\delta_{\Z_3}+6\delta_{\Z_4}+
24\delta_{\Z_2^2} +24\delta_{S_3}+24\delta_{D_8}+
24\delta_{A_4}+24 \delta_{S_4}.
\end{equation}
All the terms in the sum can be computed as above, except 
those corresponding to $D_8=\Z_4\rtimes \Z_2$, 
and $S_4=\Z_2^2\rtimes \Z_3\rtimes \Z_2$, for which 
other techniques are needed.
\end{remark}

\subsection{Normal subgroups of finite index}
\label{subsec:normal}
For each positive integer $k$, let
$a_{k}^{\nor}(G)$ be the number of index $k$, normal
subgroups of $G$.  We then have:
\begin{equation}
\label{eq:aknorm}
\a_k^{\nor}(G)=\sum_{\abs{\Gamma}=k} \delta_{\Gamma} (G).
\end{equation}
Using this formula, and our previous formulas for the 
Hall invariants, we can compute $\a_k^{\nor}(G)$ in terms 
of homological data, provided $k$ has at most two factors. 

\begin{theorem}
\label{thm:aknorm}
Let $G$ be a finitely-presented group. 
\begin{enumerate}
\item If $p$ is prime, then
\begin{align*}
a_p^{\nor}(G) &=\tfrac{p^{n}-1}{p-1}\\
a_{p^2}^{\nor}(G) &=
\tfrac{(p^{n}-1)(p^{n-1}-1)}{(p^2-1)(p-1)}+
\tfrac{p^{n-1}(p^{m}-1)}{p-1}
\end{align*}
where $n=b_1^{(p)}(G)$ and 
$m=\dim_{\,\Z_p} (p\cdot H_1(G,\Z_{p^2}))\otimes \Z_p$. 

\item If $p$ and $q$ are distinct primes, then
\begin{equation*}
a_{pq}^{\nor}(G)=
\begin{cases}
\tfrac{(p^{n}-1)(q^{m}-1)}{(p-1)(q-1)}
\quad &\text{if}\quad p\nmid q-1
\\[3pt]
\tfrac{(p^{n}-1)(q^{m}-1)}{(p-1)(q-1)}+
\tfrac{p-1}{q-1} \sum_{d\ge 1}\b_{p,d}^{(q)}(G) (q^{d}-1) 
\quad &\text{if}\quad p\mid q-1
\end{cases}
\end{equation*}
where $n=b_1^{(p)}(G)$ and $m=b_1^{(q)}(G)$.
\end{enumerate}
\end{theorem}

\begin{proof}
The only group of order $p$ is $\Z_p$; 
the only groups of order $p^2$ are $\Z_{p}^{2}$ and $\Z_{p^2}$; 
the only groups of order $pq$ are $\Z_{pq}$ (if $p\nmid q-1$), 
and $\Z_{pq}$ and $M_{p,q}$ (if $p\mid q-1$). The formulas follow from 
\eqref{eq:aknorm} and Theorems~\ref{thm:abelcount} and \ref{thm:delmpqs}. 
\end{proof}

Note that these formulas compute $a_k^{\nor}$ for all $k\le 15$, 
except for $k=8$ and $k=12$.  To compute $a_8^{\nor}$, one would need 
to know $\delta_{D_8}$ and $\delta_{Q_8}$, where $Q_8$ is the quaternion 
group; for $a_{12}^{\nor}$, one would need  
$\delta_{D_{12}}$ and $\delta_{D'_{12}}$, where 
$D'_{12}=\Z_3\rtimes \Z_4$ is the dicylic group of order $12$. 

\begin{remark}
\label{rem:dirichlet}
We may also define $\a_{k}(G)$ to be the number of index $k$, normal
subgroups $K\nor G$, with $G/K$ abelian.  That is, 
\begin{equation}
\label{eq:alphainv}
\a_k(G)=\sum_{\genfrac{}{}{0pt}{}{\Gamma
\text{ abelian}}{\abs{\Gamma}=k}}
\delta_{\Gamma} (G).
\end{equation}
Clearly, $\a_k(G)=a_k(H_1(G))$.  In particular, if 
$H_1(G)=\Z^n$, then $\a_k(G)=a_k(\Z^n)$ is given 
by the recursion \eqref{eq:bushreiner}.
\end{remark}

Finally, let $c_{k}(G)$ be the number of conjugacy classes of index $k$
subgroups of $G$. If $p$ is a prime, then clearly 
$a_p(G) = pc_p(G) - (p-1) a_p^{\nor}(G)$. 
Hence, if $n=b_1^{(p)}(G)$, we have:
\begin{equation}
\label{eq:apcp}
    c_p(G) = \tfrac{p^n+a_p(G)-1}{p}.
\end{equation}

\begin{remark}
The following formula of Stanley \cite[(5.125)]{St} holds:
$a_k(G\times \Z) = \sum_{d\mid k} d c_k(G)$.  
Hence, if $p$ is a prime, and $n=b_1^{(p)}(G)$, we have:
\begin{equation}
    a_p(G\times \Z) =a_p(G) + p^n.
\end{equation}
\end{remark}

\section{Arrangements of complex hyperplanes}
\label{sec:complexarr}

A {\em (complex) hyperplane arrangement} is a finite collection
of codimension~$1$ affine subspaces in a complex vector space.
Let $\A=\{H_1,\dots, H_n\}$ be a central arrangement of
$n$ hyperplanes in $\C^{\ell}$.
A defining polynomial for $\A$ may be written as $f=f_1\cdots f_n$,
where $f_i$ are (distinct) linear forms.  Choose coordinates
$(z_1,\dots,z_{\ell})$  in $\C^{\ell}$ so that $H_n = \ker(z_{\ell})$.
The {\em decone} $\A^*=\mathbf{d}\A$ (corresponding to this choice) 
is the affine arrangement in $\C^{\ell-1}$ with defining polynomial 
$f^* =f(z_1,\dots,z_{\ell-1},1)$.  If 
$X(\A)=\C^{\ell}\setminus \bigcup_{H\in\A} H$ is the complement 
of $\A$, then $X(\A) \cong X(\A^*)\times\C^*$.

Let $G=G(\A)=\pi_1(X(\A))$ be the fundamental group of the complement of $\A$.  
Then $G(\A)\cong G(\A^*)\times \Z$. Let $m$ be the number 
of multiple points in a generic $2$-section of $\A^*$.  The group $G^*=G(\A^*)$ 
admits a finite presentation of the form 
\begin{equation*}
\label{eq:decpres}
G^*=\langle x_1,\dots,x_{n-1} \mid
\a_j(x_i)=x_i, \text{ for } j=1,\dots, m \text{ and } i=1,\dots, n-1 \rangle,
\end{equation*}
where $\a_1,\dots, \a_m$ are the ``braid monodromy" generators---pure braids 
on $n-1$ strings, acting on $F_{n-1}=\langle x_1,\dots,x_{n-1} \rangle$ via 
the Artin representation, see \cite{CSbm} for details and further references.  
In particular, $H_1(G)=\Z^n$.

Let $V_{d}(G,\K)$ be the characteristic varieties of the arrangement $\A$ 
(over the field $\K$).  If $\K=\C$, the following facts are known:
\begin{alphenum}
\item\label{a} 
The components of $V_{d}(G,\C)$ are subtori of the character torus 
$(\C^*)^{n}$, possibly translated by roots of unity (cf.~\cite{Ar}). 
  
\item\label{b}  
The tangent cone at $\mathbf{1}$ to $V_{d}(G,\C)$ coincides 
with the resonance variety $R_d(G,\C)$; thus, the components of 
$R_d(G,\C)$ are linear subspaces of $\C^n$ (cf.~\cite{CScv, Li2}).  
\end{alphenum}
We do not know whether \eqref{a} holds 
if $\C$ is replaced by a field $\K$ of positive characteristic.    
On the other hand, the first half of \eqref{b} can easily fail 
in that case.  We refer to \cite{CScv, Fa, Li2, LY} for 
methods of computing the (complex) characteristic 
and resonance varieties of hyperplane arrangements,
and to \cite{Su} for further details on the examples below.


\begin{example} 
\label{ex:braid}
Let $\A$ be the braid arrangement in $\C^3$, with defining polynomial 
$f=xyz(x-y)(x-z)(y-z)$.  
The fundamental group is $G=P_4$, the pure braid group on $4$ strands.

For any field $\K$, the variety $V_1(G,\K)\subset (\K^*)^6$ 
has five components, all $2$-dimen\-sional:
four `local' components, corresponding to triple points, and one
`non-local' component, corresponding to an (essential) neighborly partition
of the matroid.  The components meet only at the origin, 
$\mathbf{1}=(1,\dots,1)$.  Moreover,  $V_2(G,\K)=\{\mathbf{1}\}$.

Let $p$ be a prime, and $q=0$, or a prime distinct from $p$. 
From Theorem~\ref{thm:torscount}, we get:
\begin{equation*}
\beta_{p,0}^{(q)}(G)=(p+1)(p^4+p^2-4),\quad \beta_{p,1}^{(q)}(G)=5(p+1),
\end{equation*}
and $\beta_{p,d}^{(q)}(G)=0$, for $d>1$. 
Thus, $\beta_p^{(q)}(G)=(5(p+1))$, for all $q$.
Theorem \ref{thm:delmpqs} now gives:
\begin{equation*}
\delta_{M_{p,q^s}}(G)=5(p^2-1)/s.
\end{equation*}
\end{example}

\begin{example}  
\label{ex:diamond}
Let $\A$ be the realization of the non-Fano plane, with defining polynomial 
$f=xyz(x-y)(x-z)(y-z)(x+y-z)$.  The group $G=G^*\times \Z$ 
is given by the braid monodromy generators 
$\{A_{345},\ A_{125}^{A_{35}A_{45}},\
A_{14}^{A_{34}},\ A_{136},\ A_{246}^{A_{34}A_{36}}\}$, where 
$A_I$ denotes the full twist on the strands indexed by $I$, and $x^y=y^{-1}xy$. 

Let $\K$ be a field. 
The variety $V_1(G,\K)\subset (\K^*)^7$ has nine $2$-dimensional components:
$6$ corresponding to triple points, and $3$ corresponding to
braid sub-arrangements.   All components intersect at the origin.
If $\ch\K\ne 2$, the $3$ non-local components also intersect at a
point of order~$2$, belonging to $V_2$.
Thus, $\beta_{2}^{(q)}=(24,1)$, and $\beta_{p}^{(q)}=(9(p+1))$.
Hence:
\begin{equation*}
\delta_{D_{2q}}=q+25 \quad \text{and} \quad 
\delta_{M_{p,q^s}}=9(p^2-1)/s,\ \text{for}\ p>2.
\end{equation*}
\end{example}

\begin{figure}
\setlength{\unitlength}{0.7cm}
\subfigure{
\begin{minipage}[t]{0.3\textwidth}
\begin{picture}(3,3.5)(-0.5,0)
\multiput(1,0)(1,0){2}{\line(1,1){3}}
\multiput(4,0)(1,0){2}{\line(-1,1){3}}
\put(3,0){\line(0,1){3}}
\put(4.2,-0.5){\makebox(0,0){$1$}}
\put(5.2,-0.5){\makebox(0,0){$2$}}
\put(5.2,3.5){\makebox(0,0){$3$}}
\put(4.2,3.5){\makebox(0,0){$4$}}
\put(3,3.5){\makebox(0,0){$5$}}
\end{picture}
\end{minipage}
}
\subfigure{
\begin{minipage}[t]{0.3\textwidth}
\begin{picture}(3,3.5)(-0.5,0)
\multiput(1,0)(1,0){2}{\line(1,1){3}}
\multiput(4,0)(1,0){2}{\line(-1,1){3}}
\put(3,0){\line(0,1){3}}
\put(1,1.5){\line(1,0){4}}
\put(3,-0.5){\makebox(0,0){$1$}}
\put(4.2,-0.5){\makebox(0,0){$2$}}
\put(5.2,-0.5){\makebox(0,0){$3$}}
\put(5.6,1.5){\makebox(0,0){$4$}}
\put(5.2,3.5){\makebox(0,0){$5$}}
\put(4.2,3.5){\makebox(0,0){$6$}}
\end{picture}
\end{minipage}
}
\subfigure{
\begin{minipage}[t]{0.3\textwidth}
\begin{picture}(3,3.5)(-0.5,0)
\multiput(1,0)(1,0){2}{\line(1,1){3}}
\multiput(4,0)(1,0){2}{\line(-1,1){3}}
\multiput(2.5,0)(0.5,0){3}{\line(0,1){3}}
\put(4.2,-0.5){\makebox(0,0){$1$}}
\put(5.2,-0.5){\makebox(0,0){$2$}}
\put(5.2,3.5){\makebox(0,0){$3$}}
\put(4.2,3.5){\makebox(0,0){$4$}}
\put(3.5,3.5){\makebox(0,0){$5$}}
\put(3,3.5){\makebox(0,0){$6$}}
\put(2.5,3.5){\makebox(0,0){$7$}}
\end{picture}
\end{minipage}
}
\caption{\textsf{Decones of braid, non-Fano, and deleted
$\operatorname{B}_3$ arrangements}}
\label{fig:diamonda22}
\end{figure}


\begin{example} 
\label{ex:a22}
The arrangement defined by $f=xyz(x-y)(x-z)(y-z)(x-y-z)(x-y+z)$
is a deletion of the reflection arrangement of type $\operatorname{B}_3$.
The group $G$ is isomorphic
to $(F_4\rtimes_{\a}F_3)\times \Z$,  where $\a:F_3\to P_4$ is given by
$\a=\{A_{23}, \  A_{13}^{A_{23}}A_{24},\ A_{14}^{A_{24}} \}$.

The variety $V_1(G,\K)\subset (\K^*)^8$ has eleven $2$-dimen\-sional
components ($6$ corresponding to triple points, and $5$ to braid
sub-arrange\-ments), and one $1$-dimensional component,
\begin{align*}
C&=\{t_4-t_1=t_3-t_2=t_5-t_1^2=t_7-t_2^2=t_6+1=t_8+1=t_1t_2+1=0\} \\
  &=\{(t,-t^{-1},-t^{-1},t,t^2,-1,t^{-2},-1)\mid t\in \K^*\}.
\end{align*}
The variety $V_2(G,\K)$ has a $3$-dimen\-sional component (corresponding
to a quadruple point).  Let $q=\ch\K$.  
There are two cases to consider:  

If $q\ne 2$, then $V_2(G,\K)$ also has two
isolated points of order $2$. Moreover, $C$ does not pass through the origin,
though it meets the other non-local components at the two isolated 
points of $V_2$. It follows that $\beta^{(q)}_2=(27,9)$ 
and $\beta^{(q)}_p=(11(p+1),p^2+p+1)$, for $p$ odd. 

If $q=2$, then all the components of $V_d(G,\K)$ pass through the origin.  
Note that the Galois field $\K=\F_4$, obtained by adjoining to $\F_2$ 
the primitive $3^{\text{rd}}$ root of unity $\o=e^{2 \pi \ii/3}$, 
is sufficiently large with respect to $\Z_3$.  
The representation $\mu:G\to \Z_3=\F_4^*$ given by
$\mu=(\o, \o^2, \o^2, \o,  \o^2, 1, \o, 1)$ 
belongs to $C\subset V_1(G, \F_4)$, but does not belong to $V_1(G,\C)$.  
Thus, by Corollary~\ref{cor:homprime}, $b_1(K_\mu)=8$ and
$b_1^{(2)}(K_\mu)=10$.  
Moreover, $b_1^{(2)}(K_\l)=b_1(K_\l)$, unless $\l=\mu$ or $\bar\mu$.   
It follows that $\b_{p,d}^{(2)}=\beta_{p,d}$, 
except for $\b_{3,1}^{(2)}=\beta_{3,1}+1=45$. 

Using now Theorem \ref{thm:delmpqs}, we conclude: 
\begin{align*}
\delta_{D_{2q}}&=9(q+4),\quad
\delta_{\Alt_4}=110,\\
\delta_{M_{p,q^s}}&=((p^3-1)q^s+p^3+11p^2-12)/s,\ \text{otherwise}.
\end{align*}
\end{example}

\begin{remark} 
\label{ex:fs}
Using the formulas from Section~\ref{sec:lis}, we may compute 
the number of low-index subgroups of arrangement groups. 
For example, if $\abs{\A}=n$, then 
$a_3(G)=\tfrac{1}{2}(3^{n}-1) + 3\delta_{S_3}(G)$.  
Thus, the braid arrangement has $a_3=409$,  
the non-Fano plane has $a_3=1,177$, and the 
deleted $\operatorname{B}_3$ arrangement has $a_3=3,469$. 

The idea to use the count of index $3$ subgroups as an
invariant for hyperplane arrangement groups originates
with the (unpublished) work of M.~Falk and B.~Sturmfels.   
These authors considered a pair of non-lattice-isomorphic 
arrangements of $9$ planes in $\C^3$.  The respective groups 
are in fact isomorphic (see \cite[Example~7.5]{CSbm}). 
In each case, the variety $V_1$ has twelve $2$-dimensional components
($8$ corresponding to triple points, and $4$ to braid sub-arrange\-ments),
and $V_2$ has one $3$-dimen\-sional component
(corresponding to a quadruple point). Hence:  
$\delta_{S_3}=\tfrac{1}{2}(12(2^2-1)(3-1) +(2^3-1)(3^2-1))=64$, 
$a_3=10,033$, and $c_3=9,905$.
\end{remark}

\section{Arrangements of transverse planes in $\R^4$}
\label{sec:planes}

A {\em $2$-arrangement} in $\R^{4}$ is a finite collection
$\A=\{H_1,\dots,H_n\}$ of transverse planes through the origin
of $\R^4$.   In coordinates $(z,w)$ for $\R^4=\C^2$, a defining
polynomial for $\A$ may be written as $f=f_1\cdots f_n$, where
$f_i(z,w)=a_iz+b_i\bar{z}+c_iw+d_i\bar{w}$.  A generic section 
of $\A$ by an affine $3$-plane in $\R^4$ yields a 
configuration of $n$ skew lines in $\R^3$; conversely, coning 
such a configuration yields a $2$-arrangement in $\R^{4}$. 

The complement of the arrangement, $X(\A)=\R^{4}\setminus \bigcup_{i=1}^{n} H_i$, 
deform-retracts onto the complement of the link $L(\A)=\SP^3\cap
\bigcup_{i=1}^{n} H_i$.   The link $L(\A)$ is the closure of a pure braid
$\beta\in P_n$.   The fundamental group of the complement, 
$G(\A)=\pi_1(X(\A))$, has the structure of a semidirect product of free groups:  
$G(\A)=F_{n-1}\rtimes_{\xi^2} \Z$, where $\xi$ is a certain 
pure braid in $P_{n-1}$, determined by $\beta$, see \cite{MS1}.  
Moreover, $X(\A)$ is an Eilenberg-MacLane space $K(G,1)$. 

A $2$-arrangement $\A$ is called {\em horizontal} if it admits a
defining polynomial of the form $f=\prod_{i=1}^{n}(z+a_iw+b_i\bar{w})$,
with $a_i, b_i$ real.
 From the coefficients of $f$, one reads off a permutation
$\tau\in S_n$.  Conversely, given $\tau\in S_n$, choose real numbers
$a_1 < \dots < a_n$ and $b_{\t_{1}} < \dots < b_{\t_{n}}$.
Then the polynomial
$f=\prod_{i=1}^{n}{(z-\frac{a_i+b_i}{2}w-\frac{a_i-b_i}{2}\bar{w})}$
defines a horizontal arrangement, $\A(\tau)$, whose associated
permutation is $\t$. The  braid $\xi$ corresponding to
$\A=\A(\tau)$ can be combed as $\xi=\xi_2 \cdots \xi_{n-1}$, where
$\xi_j=\prod_{i=1}^{j-1} A_{i,j}^{e_{i,j}}$ and
$e_{i,j}=1$ if $\t_{i} > \t_{j}$, and
$e_{i,j}=0$, otherwise.

For $n\le 5$, all $2$-arrangements are horizontal.  For $n=6$, 
there are $4$ non-horizontal arrangements: $\LL$, $\MM$, and 
their mirror images.  These arrangements were introduced by 
Mazurovski\u{\i} in \cite{Maz}; further details about them 
can be found in \cite{MS1}. For $n=7$, there are $13$ non-horizontal 
arrangements, see~\cite{BM}.
  
The (complex) characteristic varieties of $2$-arrangement groups $G=G(\A)$
were studied in \cite{MS1}.   Note that $V_1(G,\C)$ is the
hypersurface in $(\C^*)^n$  defined by the Alexander polynomial 
of the link $L(\A)$;  see Penne~\cite{Pe} for another way
to compute this polynomial.    We refer to \cite{MS2} for 
more information on the resonance varieties of $2$-arrangements.

\subsection{Computations of Hall invariants}
\label{subsec:hall2arr}
We now show how to compute the distribution $\beta_{p}^{(q)}(G)$ 
of mod $q$ Betti numbers of index $p$ normal subgroups, 
and the Hall invariants $\delta_{M_{p,q^s}}(G)$, for some 
$2$-arrangement groups $G$.  

\begin{example}
\label{ex:a2134}
Let $\A=\A(2134)$ be the horizontal arrangement defined by 
the polynomial $f=zw(z-w)(z-2\bar{w})$. The group of the 
complement is $G=F_3\rtimes_{\xi^2}\Z$, where $\xi=A_{1,2}$.
The characteristic varieties $V_d=V_d(G,\K)$ are given by:
\begin{align*}
\label{eq:avz}
V_1&=\{t_4=1\} \cup \{t_4=t_2^2\}, \\
V_2&=\{t_4=1, t_2=-1\} \cup \{t_4=t_2=t_1=1\} \cup \{t_4=t_2=t_3=1\},\\
V_3&=\{(1,1,1,1)\}.
\end{align*}

Counting $2$-torsion points on $V_d(G,\K)$, we see that
$\Tors_{2,1}(G,\K)\setminus \Tors_{2,2}(G,\K)=
\{(-1, 1, -1, 1)\}$, and
$\Tors_{2,2}(G,\K)\setminus \Tors_{2,3}(G,\K)=
\{(1,-1,\pm 1,1), (-1,-1,\pm 1,1)$, 
$(1, 1, -1,1), (-1,1,1,1)\}$, 
provided $\ch\K\ne 2$.  
Hence, $\beta_{2}^{(q)}=(1,6)$, if $q\ne 2$.
For an odd prime $p$, the count of $p$-torsion points yields
$\beta_{p}^{(q)}=(2p^2+p-1,2)$, if $q\nmid 2p$.

Now suppose $q= 2$. 
A sufficiently large field for the group $\Z_p$ is 
$\F_{2^s}=\F_2(\zeta)$, where $s=\ord_p(2)$
and $\zeta=e^{2 \pi \ii/p}$.  
We have $V_2(G,\F_{2^s})=\{t_4=t_2=1\}$, and so  
$\Tors_{p,2}(G,\F_{2^s})$ also contains the points
$(\zeta^j, 1, \zeta, 1)$, for $0<j<p$.
Hence, $\b_{p}^{(2)}=(2p^2,p+1)$.  By Theorem~\ref{thm:delmpqs}:
\begin{align*}
\delta_{D_{2q}}&=6q+7, &&\\
\delta_{M_{p,2^s}}&=(p-1)(2p^2+(p+1)(2^s+1))/s,&\quad &\text{for $p>2$,}\\
\delta_{M_{p,q^s}}&=(p-1)(2p^2+p+2q^s+1)/s,&\quad &\text{for $p, q>2$}.
\end{align*}

\end{example}

\begin{figure}
\setlength{\unitlength}{1.38cm}
\subfigure{
\begin{minipage}[t]{0.4\textwidth}
\begin{picture}(3,3)(-1,0)
\put(0,0){\line(2,3){0.94}}
\put(1.05,1.575){\line(2,3){0.17}}
\put(1.3,1.95){\line(2,3){0.7}}
\put(0,0.5){\line(3,2){0.46}}
\put(0.75,1){\line(3,2){0.34}}
\put(1.23,1.32){\line(3,2){0.42}}
\put(1.8,1.7){\line(3,2){0.9}}
\put(0,2.5){\line(2,-1){0.39}}
\put(0.62,2.19){\line(2,-1){2}}
\put(0,3){\line(2,-3){2}}
\put(-0.3,0){\makebox(0,0){$1$}}
\put(-0.3,0.5){\makebox(0,0){$2$}}
\put(-0.3,2.5){\makebox(0,0){$3$}}
\put(-0.3,3){\makebox(0,0){$4$}}
\end{picture}
\end{minipage}
}
\subfigure{
\begin{minipage}[t]{0.4\textwidth}
\begin{picture}(3,3)(-1,0)
\put(0,0){\line(2,3){0.5}}
\put(0.7,1.05){\line(2,3){0.25}}
\put(1.05,1.575){\line(2,3){0.28}}
\put(1.4,2.1){\line(2,3){0.7}}
\put(0,0.5){\line(3,2){1.1}}
\put(1.23,1.32){\line(3,2){1.8}}
\put(0,1.935){\line(1,0){0.61}}
\put(0.8,1.935){\line(1,0){0.39}}
\put(1.365,1.935){\line(1,0){0.185}}
\put(1.8,1.935){\line(1,0){0.25}}
\put(2.25,1.935){\line(1,0){0.8}}
\put(0,2.5){\line(3,-1){0.36}}
\put(0.54,2.32){\line(3,-1){1.385}}
\put(2.07,1.81){\line(3,-1){0.9}}
\put(0,3){\line(2,-3){2}}
\put(-0.3,0){\makebox(0,0){$1$}}
\put(-0.3,0.5){\makebox(0,0){$2$}}
\put(-0.3,1.935){\makebox(0,0){$3$}}
\put(-0.3,2.5){\makebox(0,0){$4$}}
\put(-0.3,3){\makebox(0,0){$5$}}
\end{picture}
\end{minipage}
}
\caption{\textsf{Generic $3$-sections of $\A(2134)$ and $\A(31425)$}}
\label{fig:planes}
\end{figure}

\begin{example}
\label{ex:a31425}
Let $\A=\A(31425)$.  A defining polynomial 
is $f=z(z-w)(z-2w)(z+\frac{3}{2}w-\frac{5}{2}\bar{w})
(z-\frac{1}{2}w-\frac{5}{2}\bar{w})$. The
group is $G=F_4\rtimes_{\xi^2}\Z$, where 
$\xi=A_{1,3}A_{2,3}A_{2,4}$. If $\ch\K\ne 2$, 
the varieties $V_d=V_d(G,\K)$ are as follows:

{\footnotesize
\begin{align*}
V_1=\:&\{
\begin{tabular}[t]{@{}p{5.05in}}
$ t_{2}^2 t_{3}^2 t_{4}^2 
- t_{1}^2 t_{2}^2 t_{3}^2 t_{5} 
- t_{2}^2 t_{4}^2 t_{5} 
- t_{1}^2 t_{2}^2 t_{4}^2 t_{5}  
- t_{1}^2 t_{3}^2 t_{4}^2 t_{5} 
+ {t_{2}}^4 t_{5}^2 
+ t_{1}^2 t_{3}^2 t_{5}^2 
+ t_{2}^2 t_{4}^2 t_{5}^2 
- t_{1}^2 t_{2}^2 t_{5}^3 
- {t_{2}}^4 t_{4}^2 t_{5}
+ t_{2}^2 t_{3}^2 t_{5}^2 
+ t_{1}^2 t_{2}^2 t_{3}^2 t_{5}^2 
+ 2 t_{5} (
  t_{1}^2 t_{2}^2 t_{3} t_{4} 
+ t_{2}^3 t_{3} t_{4} 
- t_{1} t_{2}^3 t_{3} t_{4} 
- t_{1} t_{2}^2 t_{3} t_{4} 
+ t_{1} t_{2} t_{3}^2 t_{4} 
- t_{1}^2 t_{2} t_{3}^2 t_{4} 
- t_{2}^2 t_{3}^2 t_{4} 
+ t_{1} t_{2}^2 t_{3}^2 t_{4} 
+ t_{1} t_{2}^2 t_{4}^2 
- t_{2}^3 t_{4}^2 
+ t_{1} t_{2}^3 t_{4}^2 
- t_{1} t_{2} t_{3} t_{4}^2 
+ t_{1}^2 t_{2} t_{3} t_{4}^2 
+ t_{2}^2 t_{3} t_{4}^2 
- t_{1} t_{2}^2 t_{3} t_{4}^2 
+ t_{1} t_{2}^2 t_{3} t_{5}
- t_{1} t_{2} t_{3}^2 t_{5} 
- t_{1}^2 t_{2}^2 t_{3} t_{5} 
- t_{2}^3 t_{3} t_{5}
+ t_{1} t_{2}^3 t_{3} t_{5}
+ t_{1}^2 t_{2} t_{3}^2 t_{5}  
- t_{1} t_{2}^2 t_{3}^2 t_{5}  
- t_{1} t_{2}^2 t_{4} t_{5}
+ t_{1}^2 t_{2}^2 t_{4} t_{5}
+ t_{2}^3 t_{4} t_{5}
- t_{1} t_{2}^3 t_{4} t_{5}
+ t_{1} t_{2} t_{3} t_{4} t_{5} 
- t_{1}^2 t_{2} t_{3} t_{4} t_{5} 
- t_{2}^2 t_{3} t_{4} t_{5}
+ t_{1} t_{2}^2 t_{3} t_{4} t_{5}
) =0
\},$
\end{tabular}
\\
V_2=\:
& \{
\begin{tabular}[t]{@{}p{5.05in}}
$ 
t_{1} t_{2}^2 t_{3} + t_{1}^2 t_{3}^2 - t_{1} t_{2} t_{3}^2 + 
  t_{1}^2 t_{2} t_{3}^2 + t_{2}^2 t_{4} - t_{1} t_{2}^2 t_{4} + 
  t_{2}^3 t_{4} + t_{1} t_{3} t_{4} - t_{1}^2 t_{3} t_{4} - 
  t_{2} t_{3} t_{4} + 4 t_{1} t_{2} t_{3} t_{4} - 
  t_{1}^2 t_{2} t_{3} t_{4} - t_{2}^2 t_{3} t_{4} + 
  t_{1} t_{2}^2 t_{3} t_{4} - t_{1} t_{3}^2 t_{4} + 
  t_{1}^2 t_{3}^2 t_{4} + t_{1}^2 t_{2}^{-1} t_{3}^2 t_{4}  + 
  t_{2} t_{4}^2 - t_{1} t_{2} t_{4}^2 + t_{2}^2 t_{4}^2 + 
  t_{1} t_{3} t_{4}^2=
t_5-t_4+t_3+t_{1}^{-1}t_{4}+t_1^{-1}t_2t_4+t_{2}^{-1}t_{3}t_{4}=0\}  
\cup \{(t,1,1,1,1)\}  \cup \{(1,1,t,1,1)\} \cup 
\{ (1,1,1,t,1)\} \cup \{(t,t,t,t,1)\}  \cup
\{(1,t,t,1,1)\}  \cup
\{ (t,t,1,1,1)\} \cup
\{(1,t,t^2,t^2,t^2)\}\cup
\{(1,1,t,t,t^2) \} \cup\{ (t,t,t,t^2,t^2)\} ,
$
\end{tabular}  
\\
V_3=\:
&\{
(-1,-1,-1,\pm 1,1),
(1,-1,-1,\pm 1,1),
(1,1,-1,\pm 1,1), 
(-1,\pm 1,1,1,1), 
(1,\pm 1,1,1,1), 
\\[-4pt]
&\ \,
(1,1,1,-1,1) 
\},
\\
V_4=\:&\{
(1,1,1,1,1)
\}.
\end{align*}
}

We start by counting points of order $2$ on these varieties.   
Inspection  shows that $\Tors_{2,1}(G,\K)\setminus \Tors_{2,2}(G,\K)=
\{ (-1, 1, -1, \pm 1, 1), (1, -1, 1, -1, 1),  (-1, 1, 1, -1, 1)$, 
$(-1,-1, 1,-1, 1)\}$ 
and $\Tors_{2,2}(G,\K)=\Tors_{2,3}(G,\K)$. 
The set $\Tors_{2,3}(G,\K)$ consists of the 
$10$ points in $V_3\setminus \{\mathbf{1}\}$, except when 
$\K=\Z_3$, in which case it also contains $(-1, 1, -1, 1, -1)$.
Hence, $\beta_{2}^{(q)}=(5,0,10)$, except for $\b_{2}^{(3)}=(5,1,10)$. 
Therefore: 
\[
\delta_{S_3}=139\quad\text{and}\quad 
\delta_{D_{2q}}=5(2q^2+2q+3),\ \text{for $q>3$}.
\]

Similar computations hold for torsion points of order $3$. 
We easily see that $\beta_{3}^{(q)}=(60,10)$, if $q\ne 2, 5, \text{ or } 7$.  
On the other hand, consider the field $\K=\F_7$, which is sufficiently large 
with respect to $\Z_3$   
(if we identify the additive group of 
$\F_7$ with the multiplicative subgroup of $\C^*$ generated by 
$\zeta=e^{2 \pi \ii/7}$, we may view $\Z_3$ as the subgroup 
$\langle \zeta^2\rangle \subset \F_7^*$).  Now  
$\Tors_{3,1}(G,\F_7)\setminus (\Tors_{3,2}(G,\F_7)\cup \Tors_{3,1}(G,\C))$ 
consists of  
$
(\zeta^2, \zeta^4, 1, \zeta^4, 1), (\zeta^2, \zeta^4, \zeta^2, \zeta^4, 
\zeta^2),  (\zeta^2, 1, \zeta^2, \zeta^4, \zeta^2), 
(\zeta^2, 1, \zeta^4, 1, \zeta^2), (\zeta^2, 1, \zeta^4, \zeta^2, 1)
$,
together with their conjugates.   
Hence, $\b_{3}^{(7)}=(65,10)$. 
Similarly, $\b_{3}^{(2)}=(41,30)$ and $\b_{3}^{(5)}=(70,10)$. 
Therefore:
\[
\delta_{\Alt_4}=191, \
\delta_{M_{3,5^2}}=330, \
\delta_{M_{3,7}}=290,  \
\text{and}\ \delta_{M_{3,q^s}}=20(q^s+7)/s,\ \text{for $q>7$}.
\]
\end{example}

\begin{remark}
\label{rem:tcone}
The resonance varieties $R_d(G,\C)$ of the arrangement $\A(31425)$ 
were computed in \cite[Example~6.5]{MS2}.  Comparing the answer given there  
with the one from Example~\ref{ex:a31425}, 
we see that the variety $R_2(G,\C)$ has $10$ irreducible components, 
whereas $V_2(G,\C)$ has only $9$ components passing through the origin 
(the tenth component, which does {\em not} contain $\mathbf{1}$, 
is {\em not} a translated torus).  Thus, the tangent cone at 
$\mathbf{1}$ to $V_2(G,\C)$ is {\em strictly} contained in $R_2(G,\C)$. 
Another example where such a strict inclusion occurs (with $G$ the 
group of a certain $4$-component link) was given in \cite[\S2.3]{Mt}. 
\end{remark}

\subsection{Classification of $2$-arrangement groups}
\label{subsec:classif2arr}
The rigid isotopy classification of configurations of $n\le 7$ skew lines 
in $\R^3$ (and, thereby, of $2$-arrangements of $n\le 7$ planes in $\R^4$) 
was established by Viro~\cite{V}, Mazurovski\u{\i}~\cite{Maz}, and 
Borobia and Mazurovski\u{\i}~\cite{BM}. Clearly, if $\A$ is rigidly 
isotopic to $\A'$, or to its mirror image, then $G(\A)\cong G(\A')$. 
The converse  was established in \cite{MS1}, for $n\le 6$, using 
certain invariants derived from the characteristic varieties $V_d(G,\C)$ 
to distinguish the homotopy types of the complements.  

We now recover the homotopy-type classification from \cite{MS1}, 
extending it from $2$-arrangements of at most $6$ planes 
to horizontal arrangements of $7$ planes, by means of a pair 
of suitably chosen metabelian Hall invariants.

\begin{theorem}
\label{thm:upto7}
For the class of $2$-arrangements of $n\le 7$ planes in $\R^4$ 
(horizontal if $n=7$), the rigid isotopy-type classification,  
up to mirror images, coincides with the isomorphism-type 
classification of the fundamental groups.  
For $n\le 6$, the groups are  classified by the Hall 
invariant $\delta_{\Alt_4}$; for $n=7$, the Hall invariant 
$\delta_{S_3}$ is also needed. 
\end{theorem}

\begin{proof}
The proof is based on the computations displayed in Table~\ref{tab:hall2arr}. 
The first column lists the rigid isotopy classes (mirror pairs identified) of
$2$-arrangements $\A$ (with $n=\abs{\A}\le 6$, or $n=7$ and $\A$ horizontal), 
according to the classification by Viro, Mazurovski\u{\i}, and  Borobia 
\cite{V, Maz, BM}. 
The next two columns list the Hall invariants $\delta_{S_3}$ and 
$\delta_{\Alt_4}$ for the corresponding $2$-arrangement groups, $G=G(\A)$. 
These invariants are computed from the  varieties $V_d(G,\F_3)$ and   
$V_d(G,\F_4)$, using Theorem~\ref{thm:delmpqs}, as in
Examples~\ref{ex:a2134} and
\ref{ex:a31425}. 
\end{proof}

Note that the classification can also be achieved 
by the Hall invariant $\delta_{M_{3,7}}$ (given in the 
last column of Table~\ref{tab:hall2arr}), either singly (for $n\le 6$), 
or together with $\delta_{S_3}$ or $\delta_{\Alt_4}$ (for $n=7$ 
and $\A$ horizontal). 
Nevertheless, no combination of  these $3$ invariants is enough 
to classify the groups of non-horizontal arrangements of $n=7$ planes. 
It would be interesting to know whether other Hall invariants can 
distinguish those $13$ groups.

\begin{table}
\begin{minipage}{.49\columnwidth}
\[
\small{
\begin{array}
{|c||r|r|r|}
\hline
\A & \delta_{S_3} & \delta_{A_4} & \delta_{M_{3,7}} \\
\hline\hline
\A(123)   & 3 & 4 & 8\\
\hline
\A(1234)  & 28  & 65 & 208\\
\A(2134)  & 25  & 38 & 72\\
\hline
\A(12345) & 195  & 840 & 4,560\\
\A(21345) & 168  & 435 & 1,184\\
\A(21435) & 150  & 273 & 632\\
\A(31425) & 139  & 191 & 290\\
\hline
\A(123456) & 1,240  & 10,285 & 96,800\\
\A(213456) & 1,051  &  5,182 & 23,560\\
\A(321456) &   997  &  4,210 & 12,640\\
\A(215436) &   889  &  2,752 & 9,940\\
\A(214356) &   907  &  2,752 & 7,288\\
\A(312546) &   799  &  1,780 & 5,008\\
\A(341256) &   799  &  2,023 & 5,200\\
\A(314256) &   750  &  1,474 & 3,688\\
\A(241536) &   704  &  1,152 & 2,368\\
\LL        &   769  &  1,631 & 4,288\\
\MM        &   685  &  1,126 & 1,900\\
\hline
\A(1234567) & 7,623 & 124,124 &2,039,128 \\
\A(2134567) & 6,408 &  62,159 &  488,568 \\
\A(3214567) & 5,922 &  47,579 &  210,024 \\
\hline
\end{array}
}
\]
\end{minipage}
\begin{minipage}{.49\columnwidth}
\renewcommand{\arraystretch}{1.105}
\[
\small{
\begin{array}
{|c||r|r|r|}
\hline
\A & \delta_{S_3} & \delta_{A_4} & \delta_{M_{3,7}} \\
\hline\hline
\A(2143567) & 5,436 &  31,541 &  125,760 \\
\A(2154367) & 5,112 &  25,709 &   83,928 \\
\A(2165437) & 5,274 &  31,541 &  194,976 \\
\A(3216547) & 4,950 &  25,709 &  145,080 \\
\A(2143657) & 4,680 &  16,961 &   49,872 \\
\A(3412567) & 4,464 &  20,606 &   74,640 \\
\A(3125467) & 4,572 &  16,961 &   59,952 \\
\A(4123657) & 4,464 &  16,961 &   72,720 \\
\A(3126457) & 4,032 &  11,129 &   39,432 \\
\A(3254167) & 4,032 &  12,587 &   41,160 \\
\A(3142567) & 4,237 &  15,227 &   66,330 \\
\A(3142657) & 3,796 &   9,197 &   25,974 \\
\A(3145267) & 3,931 &  11,651 &   43,298 \\
\A(3415267) & 3,796 &  10,175 &   34,410 \\
\A(3154267) & 3,850 &  10,751 &   34,410 \\
\A(2415367) & 3,727 &   9,349 &   26,604 \\
\A(2415637) & 3,619 &   8,709 &   26,532 \\
\A(2516347) & 3,484 &   7,459 &   20,452 \\
\A(3625147) & 3,245 &   6,349 &   15,736 \\
\A(4136257) & 3,329 &   6,189 &   15,082 \\
\A(5264137) & 3,417 &   6,819 &   15,184 \\
\hline
\end{array}
}
\]
\end{minipage}
\\[8pt]
\caption{\textsf{Hall invariants of groups of $2$-arrangements in $\R^4$.}}
\label{tab:hall2arr}
\end{table}



\begin{thebibliography}{AAA}

\bibitem{Ar} D.~Arapura,
{\em  Geometry of cohomology support loci for local systems {\rm I}},
J. Alg. Geom. \textbf{6} (1997), 563--597.

\bibitem{BM} A.~Borobia, V.~Mazurovski\u{\i}, 
{\em Nonsingular configurations of $7$ lines of $\RP^3$}, 
J. Knot Theory Ramifications \textbf{6} (1997), 751--783.

\bibitem{Brd} G.~Bredon, 
{\em Sheaf theory}, Second edition, Grad. Texts in Math., vol.~170, 
Springer-Verlag, New York, 1997.

\bibitem{Br82} K.~S.~Brown,
{\em Cohomology of groups},
Corrected reprint of the 1982 original, Grad. Texts in Math., vol.~87. 
Springer-Verlag, New York, 1994. 

\bibitem{Br} \bysame,
{\em The coset poset and probabilistic zeta function of a finite group},
J. Algebra \textbf{225} (2000), 989--1012.

\bibitem{BuR} C.~Bushnell, I.~Reiner,
{\em Zeta functions of arithmetic orders and Solomon's conjectures},
Math. Z. \textbf{173} (1980), 135--161.

\bibitem{Cog} J.~Cogolludo, 
{\em Topological invariants of the complements to rational arrangements}, 
Ph.D. thesis, Univ. of Illinois at Chicago and Univ. Complutense, Madrid, 1999. 

\bibitem{CSbm} D.~Cohen, A.~Suciu,
{\em The braid monodromy of plane algebraic curves and
hyperplane arrangements}, Comment. Math. Helvetici
\textbf{72} (1997), 285--315.

\bibitem{CScv} \bysame,
{\em Characteristic varieties of arrangements},
Math. Proc. Cambridge Phil. Soc. \textbf{127} (1999), 33--53.

\bibitem{CR} C.~Curtis, I.~Reiner,
{\em Methods of representation theory. Vol.~1} (reprint
of the 1981 edition),  John Wiley \& Sons, New York, 1990.

\bibitem{Ei}  D.~Eisenbud,
{\em Commutative algebra with a view towards algebraic geometry},
Grad. Texts in Math., vol.~150, Springer-Verlag, New~York, 1995.

\bibitem{Fa} M.~Falk,
{\em Arrangements and cohomology},
Ann. Combin. \textbf{1} (1997), 135--157.

\bibitem{Fx1} R.~H.~Fox,
{\em Free differential calculus. III. Subgroups},
Ann. of Math. \textbf{64} (1956), 407--419.

\bibitem{Fx2} \bysame,
{\em Metacyclic invariants of knots and links},
Canad. J. Math. \textbf{22} (1970), 193--201.

\bibitem{gap} The GAP~Group, \emph{{GAP--Groups, Algorithms,
and Programming, Version~4.1}},  Aachen, St.~Andrews, 1999;
available at \texttt{\href{http://www-gap.dcs.st-and.ac.uk/~gap}
{http://www-gap.dcs.st-and.ac.uk/\~{}gap.}}

\bibitem{GS}  D.~Grayson, M.~Stillman,
{\em Macaulay~$2$: a software system for algebraic
geometry and commutative algebra};
available at \texttt{\href{http://www.math.uiuc.edu/Macaulay2}
{http://www.math.uiuc.edu/Macaulay2.}}

\bibitem{Ha} M.~Hall,
{\em Subgroups of finite index in free groups},
Canad. J. Math \textbf{1} (1949), 187--190.

\bibitem{HaP} P.~Hall,
{\em The Eulerian functions of a group},
Quart. J. Math \textbf{7} (1936), 134--151.

\bibitem{He1}  J.~Hempel,
{\em Homology of coverings}, Pacific J. Math. \textbf{112} (1984),
83--113.

\bibitem{He2} \bysame,
{\em Homology of branched coverings of $3$-manifolds},
Canad. J. Math. \textbf{44} (1992), 119--134.

\bibitem{Hr} E.~Hironaka,
{\em Alexander stratifications of character varieties},
Annales de l'Institut Fourier (Grenoble) \textbf{47} (1997), 555--583.

\bibitem{KCL} J.~H.~Kwak, J.-H~Chun, J.~Lee, 
{\em  Enumeration of regular graph coverings having finite
abelian covering transformation groups}, SIAM J.
Discrete Math. \textbf{11} (1998), 273--285. 

\bibitem{Li1} A.~Libgober,
{\em On the homology of finite abelian coverings},
Topology. Appl. \textbf{43} (1992), 157--166.

\bibitem{Li2} \bysame,
{\em Characteristic varieties of algebraic curves}, 
\texttt{\href{http://front.math.ucdavis.edu/math.AG/9801070}%
{arXiv:math.AG/9801070.}}

\bibitem{Li3}  \bysame, 
{\em First order deformations for rank one local systems 
with non vanishing cohomology}, Topology Appl. (to appear).

\bibitem{LY}  A.~Libgober, S.~Yuzvinsky,
{\em Cohomology of the Orlik-Solomon algebras and local systems},
Compositio Math. \textbf{21} (2000), 337--361. 

\bibitem{Lin}  D.~Lind, 
{\em A zeta function for $\Z^d$-actions}, 
In: Ergodic theory of $\Z^d$ actions (Warwick, 1993--1994), 
London Math. Soc. Lecture Note Ser., vol.~228, 
Cambridge Univ. Press, Cambridge, 1996, pp.~433--450.

\bibitem{Liv} C.~Livingston, 
{\em Lifting representations of knot groups},
J. Knot Theory Ramifications \textbf{4} (1995), 225--234. 

\bibitem{Lu}  A.~Lubotzky,
{\em Counting finite index subgroups},
In: Groups '93, Galway/St.~Andrews, Vol.~2,
London Math. Soc. Lecture Note Ser., vol.~212,
Cambridge Univ. Press, Cambridge, 1995, pp.~368--404.

\bibitem{Mac} I.~G.~Macdonald,
{\em Symmetric functions and Hall polynomials},
Second edition. With contributions by A. Zelevinsky, 
Oxford Math. Monographs, Oxford Univ. Press, New York, 1995.

\bibitem{Mt} D.~Matei,
{\em Fundamental groups of links and arrangements:
characteristic varieties, resonance varieties, and finite index subgroups},
Ph.D. thesis, Northeastern Univ., Boston, MA, 1999.

\bibitem{MS1} D.~Matei, A.~Suciu,
{\em Homotopy types of complements of $2$-arrangements in $\R^4$},
Topology \textbf{39} (2000), 61--88.

\bibitem{MS2} \bysame,
{\em Cohomology rings and nilpotent quotients of real and
complex arrangements},
In: Arrange\-ments--Tokyo 1998 (M.~Falk, H.~Terao, eds.),
Adv. Stud. Pure Math., vol.~27,  Math. Soc. Japan, Kinokuniya, 
Tokyo, 2000, pp.~185--215.

\bibitem{Maz} V.~Mazurovski\u{\i},
{\em Configurations of six skew lines} (Russian),
Zap. Nauchn. Sem. Leningrad Otdel. Mat. Inst. Steklov (LOMI)
\textbf{167} (1988), 121--134; translation in J.~Soviet Math.
\textbf{52} (1990), 2825--2832.

\bibitem{Me} A.~D.~Mednykh,
{\em Unramified coverings of compact Riemann surfaces} (Russian),
Dokl. Akad. Nauk SSSR \textbf{244} (1979), 529--532;
English translation: Soviet Math. Dokl. \textbf{20} (1979), 85--88.

\bibitem{Pe} R.~Penne, 
{\em The Alexander polynomial of a configuration of skew lines in $3$-space},  
Pacific J. Math. \textbf{186} (1998), 315--348. 

\bibitem{Rob} D.~Robinson,
{\em A course in the theory of groups}, Grad. Texts in Math., 
vol.~80, Springer-Verlag, New York, 1996.

\bibitem{Rom} S.~Roman,
{\em Field theory}, Grad. Texts in Math., vol.~158,
Springer-Verlag, New York, 1995.

\bibitem{Sa} M.~Sakuma,
{\em Homology of abelian coverings of links and spatial graphs},
Canad. J. Math. \textbf{47} (1995), 201--224.

\bibitem{St} R.~Stanley,
{\em Enumerative combinatorics. Vol.~\textup{2}}, Cambridge Studies in
Advanced Math., vol.~62, Cambridge Univ. Press, Cambridge, 1999.

\bibitem{Su} A.~Suciu,
{\em Translated tori in the characteristic varieties of complex
hyperplane arrangements}, Topology. Appl. (to appear),   
\texttt{\href{http://front.math.ucdavis.edu/math.AG/9912227}
{arXiv:math.AG/9912227.}}

\bibitem{V} O.~Viro, 
{\em Topological problems concerning lines and points of 
three-dimensional space},
Dokl. Akad. Nauk. SSSR \textbf{284} (1985), 1049--1052; English transl., 
Soviet Math. Dokl. \textbf{32} (1985), 528--531. 

\bibitem{We} L.~Weisner,
{\em Some properties of prime-power groups},
Trans. Amer. Math. Soc. \textbf{38} (1935), 485--492.

\end{thebibliography}
\end{document}